# Ranking Median Regression:
# Learning to Order through Local Consensus

**Stephan Clémençon**  STEPHAN.CLEMENCON@TELECOM-PARISTECH.FR
**Anna Korba**  ANNA.KORBA@TELECOM-PARISTECH.FR
*LTCI, Télécom ParisTech, Université Paris-Saclay*
*Paris, France*

**Eric Sibony**  ERIC.SIBONY@GMAIL.COM
*Shift Technologies*
*Paris, France*

## Abstract

This article is devoted to the problem of predicting the value taken by a random permutation $\Sigma$, describing the preferences of an individual over a set of numbered items $\{1, \ldots, n\}$ say, based on the observation of an input/explanatory r.v. $X$ (*e.g.* characteristics of the individual), when error is measured by the Kendall $\tau$ distance. In the probabilistic formulation of the 'Learning to Order' problem we propose, which extends the framework for statistical Kemeny ranking aggregation developed in Korba et al. (2017), this boils down to recovering conditional Kemeny medians of $\Sigma$ given $X$ from i.i.d. training examples $(X_1, \Sigma_1), \ldots, (X_N, \Sigma_N)$. For this reason, this statistical learning problem is referred to as *ranking median regression* here. Our contribution is twofold. We first propose a probabilistic theory of ranking median regression: the set of optimal elements is characterized, the performance of empirical risk minimizers is investigated in this context and situations where fast learning rates can be achieved are also exhibited. Next we introduce the concept of local consensus/median, in order to derive efficient methods for ranking median regression. The major advantage of this local learning approach lies in its close connection with the widely studied Kemeny aggregation problem. From an algorithmic perspective, this permits to build predictive rules for ranking median regression by implementing efficient techniques for (approximate) Kemeny median computations at a local level in a tractable manner. In particular, versions of $k$-nearest neighbor and tree-based methods, tailored to ranking median regression, are investigated. Accuracy of piecewise constant ranking median regression rules is studied under a specific smoothness assumption for $\Sigma$'s conditional distribution given $X$. The results of various numerical experiments are also displayed for illustration purpose.

**Keywords:** Consensus ranking, empirical risk minimization, fast rates, Kemeny median, local learning, nearest-neighbors, decision trees, predictive learning, ranking aggregation

## 1. Introduction

The machine-learning problem considered in this paper is easy to state. Given a vector $X$ of attributes describing the characteristics of an individual, the goal is to predict her preferences over a set of $n \geq 1$ numbered items, indexed by $\{1, \ldots, n\}$ say, modelled as a random permutation $\Sigma$ in $\mathfrak{S}_n$. Based on the observation of independent copies of the random pair $(X, \Sigma)$, the task consists in building a predictive function $s$ that maps any point $X$ in the input space to a permutation $s(X)$, the accuracy of the prediction being measured by means of a certain distance between $\Sigma$ and $s(X)$, the Kendall $\tau$ distance typically. This problem is of growing importance these days, since users with declared characteristics express their preferences through more and more devices/interfaces (*e.g.* social surveys, web activities...). This article proposes a probabilistic analysis of this statistical learning problem: optimal predictive rules are exhibited and (fast) learning rate bounds for empirical risk minimizers are established in particular. However, truth should be said, this problem is more difficult to solve in practice than other supervised learning problems such as classification or regression, due to the (structured) nature of the output space. The symmetric group is not a vector space and its elements cannot be defined by means of simple operations, such as thresholding some real valued function, like in classification. Hence, it is far from straightforward in general to find analogues of methods for distribution-free regression or classification consisting in expanding the decision function using basis functions in a flexible dictionary (*e.g.* splines, wavelets) and fitting its coefficients from training data, with the remarkable exception of techniques building *piecewise constant* predictive functions, such as the popular nearest-neighbor method or the celebrated CART algorithm, see Breiman et al. (1984). Indeed, observing that, when $X$ and $\Sigma$ are independent, the best predictions for $\Sigma$ are its Kemeny medians (*i.e.* any permutation that is closest to $\Sigma$ in expectation, see the probabilistic formulation of ranking aggregation in Korba et al. (2017)), we consider *local learning* approaches in this paper. Conditional Kemeny medians of $\Sigma$ at a given point $X = x$ are relaxed to Kemeny medians within a region $\mathcal{C}$ of the input space containing $x$ (*i.e.* local consensus), which can be computed by applying locally any ranking aggregation technique. Beyond computational tractability, it is motivated by the fact that, as shall be proved in this paper, the optimal ranking median regression rule can be well approximated by piecewise constants under the hypothesis that the pairwise conditional probabilities $\mathbb{P}\{\Sigma(i) < \Sigma(j) \mid X = x\}$, with $1 \leq i < j \leq n$, are Lipschitz. Two methods based on the notion of *local Kemeny consensus* are investigated here. The first technique is a version of the popular nearest neighbor method tailored to ranking median regression, while the second one, refered to as the CRIT algorithm (standing for 'Consensus RankIng Tree'), produces, by successive data-driven refinements, an adaptive partitioning of the input space $\mathcal{X}$ formed of regions, where the $\Sigma_i$'s exhibit low variability. Like CART, the recursive learning process CRIT can be described by a binary tree, whose terminal leafs are associated with the final regions. It can be seen as a variant of the methodology introduced in Yu et al. (2010): we show here that the node impurity measure they originally propose can be related to the local ranking median regression risk, the sole major difference being the specific computationally effective method we consider for computing local predictions, *i.e.* for assigning permutations to terminal nodes. Beyond approximation theoretic arguments, its computational feasability and the advantages of the predictive rules it produces regarding interpretability or aggregation are also discussed to support the use of piecewise constants. The results of various numerical experiments are also displayed in order to illustrate the approach we propose.

The paper is organized as follows. In section 2, concepts related to statistical (Kemeny) ranking aggregation are briefly recalled, the ranking predictive problem being next formulated as an extension of the latter and studied from a theoretical perspective. A probabilistic theory of ranking median regression is developed in section 3. In section 4, approximation of optimal predictive rules by piecewise constants is investigated as well as two local learning methods for solving ranking median regression. The results of illustrative numerical experiments are presented in section 5. Technical proofs and further details can be found in the Appendix.

## 2. Preliminaries

As a first go, we start with recalling the metric approach to consensus ranking and give next a rigorous statistical formulation of ranking aggregation. We also establish statistical guarantees for the generalization capacity of an alternative to the empirical Kemeny median technique studied in Korba et al. (2017), which can be much more easily computed in certain situations and which the algorithms we propose in the subsequent section highly rely on. From the angle embraced in this paper, *ranking median regression* is then viewed as an extension of statistical ranking aggregation. Here and throughout, the indicator function of any event $\mathcal{E}$ is denoted by $\mathbb{I}\{\mathcal{E}\}$, the Dirac mass at any point $a$ by $\delta_a$ and the cardinality of any finite set $E$ by $\#E$. Let $n \geq 1$, the set of permutations of $[\![n]\!] = \{1, \ldots, n\}$ is denoted by $\mathfrak{S}_n$.

### 2.1 Consensus Ranking: Probabilistic Framework and Statistical Setup

Throughout the article, a ranking on a set of items indexed by $[\![n]\!]$ is seen as the permutation $\sigma \in \mathfrak{S}_n$ that maps any item $i$ to its rank $\sigma(i)$. Given a collection of $N \geq 1$ rankings $\sigma_1, \ldots, \sigma_N$, the goal of consensus ranking, also referred to as ranking aggregation sometimes, is to find $\sigma^* \in \mathfrak{S}_n$ that best summarizes it. A popular way of tackling this problem, the metric-based consensus approach, consists in solving:

$$\min_{\sigma \in \mathfrak{S}_n} \sum_{i=1}^{N} d(\sigma, \sigma_i), \tag{1}$$

where $d(.,.)$ is a certain metric on $\mathfrak{S}_n$. As the set $\mathfrak{S}_n$ is of finite cardinality, though not necessarily unique, such a barycentric permutation, called *consensus/median ranking*, always exists. In Kemeny ranking aggregation, the version of this problem the most widely documented in the literature, one considers the number of pairwise disagreements as metric, namely the Kendall's $\tau$ distance, see Kemeny (1959): $\forall (\sigma, \sigma') \in \mathfrak{S}_n^2$,

$$d_\tau(\sigma, \sigma') = \sum_{i<j} \mathbb{I}\{(\sigma(i) - \sigma(j))(\sigma'(i) - \sigma'(j)) < 0\}. \tag{2}$$

The problem (1) can be viewed as a $M$-estimation problem in the probabilistic framework stipulating that the collection of rankings to be aggregated/summarized is composed of $N \geq 1$ independent copies $\Sigma_1, \ldots, \Sigma_N$ of a generic r.v. $\Sigma$, defined on a probability space $(\Omega, \mathcal{F}, \mathbb{P})$ and drawn from an unknown probability distribution $P$ on $\mathfrak{S}_n$ (*i.e.* $P(\sigma) = \mathbb{P}\{\Sigma = \sigma\}$ for any $\sigma \in \mathfrak{S}_n$). Just like a median of a real valued r.v. $Z$ is any scalar closest to $Z$ in the $L_1$ sense, a (true) median of distribution $P$ w.r.t. a certain metric $d$ on $\mathfrak{S}_n$ is any solution of the minimization problem:

$$\min_{\sigma \in \mathfrak{S}_n} L_P(\sigma), \tag{3}$$

where $L_P(\sigma) = \mathbb{E}_{\Sigma \sim P}[d(\Sigma, \sigma)]$ denotes the expected distance between any permutation $\sigma$ and $\Sigma$. In this framework, statistical ranking aggregation consists in recovering a solution $\sigma^*$ of this minimization problem, plus an estimate of this minimum $L_P^* = L_P(\sigma^*)$, as accurate as possible, based on the observations $\Sigma_1$, ..., $\Sigma_N$. A median permutation $\sigma^*$ can be interpreted as a central value for distribution $P$, while the quantity $L_P^*$ may be viewed as a dispersion measure. Like problem (1), the minimization problem (3) has always a solution but can be multimodal. However, the functional $L_P(.)$ is unknown in practice, just like distribution $P$. Suppose that we would like to avoid rigid parametric assumptions on $P$ such as those stipulated by the Mallows model, see Mallows (1957), and only have access to the dataset $\{\Sigma_1, \ldots, \Sigma_N\}$ to find a reasonable approximant of a median. Following the Empirical Risk Minimization (ERM) paradigm (see *e.g.* Vapnik, 2000), one substitutes in (3) the quantity $L(\sigma)$ with its statistical version

$$\widehat{L}_N(\sigma) = \frac{1}{N} \sum_{i=1}^N d(\Sigma_i, \sigma) = L_{\widehat{P}_N}(\sigma), \tag{4}$$

where $\widehat{P}_N = (1/N) \sum_{i=1}^N \delta_{\Sigma_i}$ denotes the empirical measure. The performance of empirical consensus rules, namely solutions $\widehat{\sigma}_N$ of $\min_{\sigma \in \mathfrak{S}_n} \widehat{L}_N(\sigma)$, has been investigated in Korba et al. (2017). Precisely, rate bounds of order $O_\mathbb{P}(1/\sqrt{N})$ for the excess of risk $L_P(\widehat{\sigma}_N) - L_P^*$ in probability/expectation have been established and proved to be sharp in the minimax sense, when $d$ is the Kendall's $\tau$ distance. Whereas problem (1) is NP-hard in general (see Hudry (2008) for instance), in the Kendall's $\tau$ case, exact solutions, referred to as *Kemeny medians*, can be explicited when the pairwise probabilities $p_{i,j} = \mathbb{P}\{\Sigma(i) < \Sigma(j)\}$, $1 \leq i \neq j \leq n$, fulfill the following property, referred to as *stochastic transitivity*.

**Definition 1** *The probability distribution $P$ on $\mathfrak{S}_n$ is stochastically transitive iff*

$$\forall (i,j,k) \in [\![n]\!]^3: \ p_{i,j} \geq 1/2 \ \text{and} \ p_{j,k} \geq 1/2 \ \Rightarrow \ p_{i,k} \geq 1/2.$$

*If, in addition, $p_{i,j} \neq 1/2$ for all $i < j$, $P$ is said to be strictly stochastically transitive.*

When stochastic transitivity holds true, the set of Kemeny medians (see Theorem 5 in Korba et al. (2017)) is the (non empty) set

$$\{\sigma \in \mathfrak{S}_n: \ (p_{i,j} - 1/2)(\sigma(j) - \sigma(i)) > 0 \text{ for all } i < j \text{ s.t. } p_{i,j} \neq 1/2\}, \tag{5}$$

the minimum is given by

$$L_P^* = \sum_{i<j} \min\{p_{i,j}, 1 - p_{i,j}\} = \sum_{i<j} \{1/2 - |p_{i,j} - 1/2|\} \tag{6}$$

and, for any $\sigma \in \mathfrak{S}_n$, $L_P(\sigma) - L_P^* = 2 \sum_{i<j} |p_{i,j} - 1/2| \cdot \mathbb{I}\{(\sigma(i) - \sigma(j))(p_{i,j} - 1/2) < 0\}$. If a strict version of stochastic transitivity is fulfilled, we denote by $\sigma_P^*$ the Kemeny median which is unique and given by the Copeland ranking:

$$\sigma_P^*(i) = 1 + \sum_{k \neq i} \mathbb{I}\{p_{i,k} < 1/2\} \text{ for } 1 \leq i \leq n. \tag{7}$$

Recall also that examples of stochastically transitive distributions on $\mathfrak{S}_n$ are numerous and include most popular parametric models such as Mallows or Bradley-Terry-Luce-Plackett models, see *e.g.* Mallows (1957) or Plackett (1975).

**Remark 2** (MEASURING DISPERSION) *As noticed in Korba et al. (2017), an alternative measure of dispersion is given by $\gamma(P) = (1/2)\mathbb{E}[d(\Sigma, \Sigma')]$, where $\Sigma'$ is an independent copy of $\Sigma$. When $P$ is not strictly stochastically transitive, the latter can be much more easily estimated than $L_P^*$, insofar as no (approximate) median computation is needed: indeed, a natural estimator is given by the U-statistic $\widehat{\gamma}_N = 2/(N(N-1)) \sum_{i<j} d(\Sigma_i, \Sigma_j)$. For this reason, this empirical dispersion measure will be used as a splitting criterion in the partitioning algorithm proposed in subsection 4.3. Observe in addition that $\gamma(P) \leq L_P^* \leq 2\gamma(P)$ and, when $d = d_\tau$, we have $\gamma(P) = \sum_{i<j} p_{i,j}(1 - p_{i,j})$.*

We denote by $\mathcal{T}$ the set of strictly stochastically transitive distributions on $\mathfrak{S}_n$. Assume that the underlying distribution $P$ belongs to $\mathcal{T}$ and verifies a certain low-noise condition **NA**$(h)$, defined for $h > 0$ by:

$$\min_{i<j} |p_{i,j} - 1/2| \geq h \qquad (8)$$

It is shown in Korba et al. (2017) that the empirical distribution $\widehat{P}_N$ is strictly stochastically transitive as well, with overwhelming probability, and that the expectation of the excess of risk of empirical Kemeny medians decays at an exponential rate, see Proposition 14 therein. In this case, the nearly optimal solution $\sigma_{\widehat{P}_N}^*$ can be made explicit and straightforwardly computed using Eq. (7) based on the empirical pairwise probabilities

$$\widehat{p}_{i,j} = \frac{1}{N} \sum_{k=1}^{N} \mathbb{I}\{\Sigma_k(i) < \Sigma_k(j)\}, \; i < j.$$

Otherwise, solving the NP-hard problem $\min_{\sigma \in \mathfrak{S}_n} L_{\widehat{P}_N}(\sigma)$ requires to get an empirical Kemeny median. However, as can be seen by examining the argument of Proposition 14's proof in Korba et al. (2017), the exponential rate bound holds true for any candidate $\widetilde{\sigma}_N$ in $\mathfrak{S}_n$ that coincides with $\sigma_{\widehat{P}_N}^*$ when the empirical distribution lies in $\mathcal{T}$.

**Theorem 3** *(Korba et al. (2017), Proposition 14) Suppose that $P \in \mathcal{T}$ and fulfills condition* **NA**$(h)$. *On the event $\{\widehat{P}_N \in \mathcal{T}\}$, define $\widetilde{\sigma}_{\widehat{P}_N} = \sigma_{\widehat{P}_N}^*$ and set $\widetilde{\sigma}_{\widehat{P}_N} = \sigma$, for $\sigma \in \mathfrak{S}_n$ arbitrarily chosen, on the complementary event. Then, for all $\delta \in (0,1)$, we have with probability at least $1 - \delta$: $\forall N \geq 1$,*

$$L_P(\widetilde{\sigma}_{\widehat{P}_N}) - L_P^* \leq \frac{n^2(n-1)^2}{8} \exp\left(-\frac{N}{2} \log\left(\frac{1}{1-4h^2}\right)\right).$$

In practice, when $\widehat{P}_N$ does not belong to $\mathcal{T}$, we propose to consider as a pseudo-empirical median any permutation $\widetilde{\sigma}_{\widehat{P}_N}^*$ that ranks the objects as the empirical Borda count:

$$\left(\sum_{k=1}^{N} \Sigma_k(i) - \sum_{k=1}^{N} \Sigma_k(j)\right) \cdot \left(\widetilde{\sigma}_{\widehat{P}_N}^*(i) - \widetilde{\sigma}_{\widehat{P}_N}^*(j)\right) > 0 \text{ for all } i < j \text{ s.t. } \sum_{k=1}^{N} \Sigma_k(i) \neq \sum_{k=1}^{N} \Sigma_k(j),$$

breaking possible ties in an arbitrary fashion. Alternative choices could also be guided by least-squares approximation of the empirical pairwise probabilities, as revealed by the following analysis.

## 2.2 Best strictly stochastically transitive approximation

We suppose that $P \in \mathcal{T}$. If the empirical estimation $\widehat{P}_N$ of $P$ does not belong to $\mathcal{T}$, a natural strategy would consist in approximating it by a strictly stochastically transitive probability distribution $\widetilde{P}$ as accurately as possible (in a sense that is specified below) and consider the (unique) Kemeny median of the latter as an approximate median for $\widehat{P}_N$ (for $P$, respectively). It is legitimated by the result below, whose proof is given in the Appendix.

**Lemma 4** *Let $P'$ and $P''$ be two probability distributions on $\mathfrak{S}_n$.*

(i) *Let $\sigma_{P''}$ ba any Kemeny median of distribution $P''$. Then, we have:*

$$L_{P'}^* \leq L_{P'}(\sigma_{P''}) \leq L_{P'}^* + 2\sum_{i<j} |p'_{i,j} - p''_{i,j}|, \qquad (9)$$

*where $p'_{i,j} = \mathbb{P}_{\Sigma \sim P'}\{\Sigma(i) < \Sigma(j)\}$ and $p''_{i,j} = \mathbb{P}_{\Sigma \sim P''}\{\Sigma(i) < \Sigma(j)\}$ for any $i < j$.*

(ii) *Suppose that $(P', P'') \in \mathcal{T}^2$ and set $h = \min_{i<j} |p''_{i,j} - 1/2|$. Then, we have:*

$$d_\tau(\sigma_{P'}^*, \sigma_{P''}^*) \leq (1/h) \sum_{i<j} |p'_{i,j} - p''_{i,j}|. \qquad (10)$$

We go back to the approximate Kemeny aggregation problem and suppose that it is known *a priori* that the underlying probability $P$ belongs to a certain subset $\mathcal{T}'$ of $\mathcal{T}$, on which the quadratic minimization problem

$$\min_{P' \in \mathcal{T}'} \sum_{i<j} (p'_{i,j} - \widehat{p}_{i,j})^2 \qquad (11)$$

can be solved efficiently (by orthogonal projection typically, when $\mathcal{T}'$ is a vector space or a convex set, up to an appropriate reparametrization). In Jiang et al. (2010), the case

$$\mathcal{T}' = \{P': \; (p_{i,j} - 1/2) + (p_{j,k} - 1/2) + (p_{k,i} - 1/2) = 0 \text{ for all 3-tuple } (i,j,k)\} \subset \mathcal{T}$$

has been investigated at length in particular. Denoting by $\widetilde{P}$ the solution of (11), we deduce from Lemma 4 combined with Cauchy-Schwarz inequality that

$$L_{\widehat{P}_N}^* \leq L_{\widehat{P}_N}(\sigma_{\widetilde{P}}^*) \leq L_{\widehat{P}_N}^* + \sqrt{2n(n-1)} \left( \sum_{i<j} (\widetilde{p}_{i,j} - \widehat{p}_{i,j})^2 \right)^{1/2}$$

$$\leq L_{\widehat{P}_N}^* + \sqrt{2n(n-1)} \left( \sum_{i<j} (p_{i,j} - \widehat{p}_{i,j})^2 \right)^{1/2},$$

where the final upper bound can be easily shown to be of order $O_\mathbb{P}(1/\sqrt{N})$.

## 2.3 Predictive Ranking and Statistical Conditional Models

We suppose now that, in addition to the ranking $\Sigma$, one observes a random vector $X$, defined on the same probability space $(\Omega, \mathcal{F}, \mathbb{P})$, valued in a feature space $\mathcal{X}$ (of possibly high dimension, typically a subset of $\mathbb{R}^d$ with $d \geq 1$) and modelling some information hopefully useful to predict $\Sigma$ (or at least to recover some of its characteristics). The joint distribution of the r.v. $(\Sigma, X)$ is described by $(\mu, P_X)$, where $\mu$ denotes $X$'s marginal distribution and $P_X$ means the conditional probability distribution of $\Sigma$ given $X$: $\forall \sigma \in \mathfrak{S}_n$, $P_X(\sigma) = \mathbb{P}\{\Sigma = \sigma \mid X\}$ almost-surely. The marginal distribution of $\Sigma$ is then $P(\sigma) = \int_{\mathcal{X}} P_x(\sigma)\mu(x)$. Whereas ranking aggregation methods (such as that analyzed in Theorem 3 from a statistical learning perspective) applied to the $\Sigma_i$'s would ignore the information carried by the $X_i$'s for prediction purpose, our goal is to learn a predictive function $s$ that maps any point $X$ in the input space to a permutation $s(X)$ in $\mathfrak{S}_n$. This problem can be seen as a generalization of multiclass classification and has been referred to as *label ranking* in Tsoumakas et al. (2009) and Vembu and Gärtner (2010) for instance. Some approaches are rule-based (see Gurrieri et al. (2012)), while certain others adapt classic algorithms such as those investigated in section 4 to this problem (see Yu et al. (2010)), but most of the methods documented in the literature rely on parametric modeling (see Cheng and Hüllermeier (2009), Cheng et al. (2009), Cheng et al. (2010)). In parallel, several authors proposed to model explicitly the dependence of the parameter $\theta$ w.r.t. the covariate $X$ and rely next on MLE or Bayesian techniques to compute a predictive rule. One may refer to Rendle et al. (2009) or Lu and Negahban (2015). In contrast, the approach we develop in the next section aims at formulating the ranking regression problem, free of any parametric assumptions, in a general statistical framework.

## 3. Ranking Median Regression

Let $d$ be a metric on $\mathfrak{S}_n$, assuming that the quantity $d(\Sigma, \sigma)$ reflects the cost of predicting a value $\sigma$ for the ranking $\Sigma$, one can formulate the predictive problem that consists in finding a measurable mapping $s : \mathcal{X} \to \mathfrak{S}_n$ with minimum prediction error:

$$\mathcal{R}(s) = \mathbb{E}_{X \sim \mu}[\mathbb{E}_{\Sigma \sim P_X}[d(s(X), \Sigma)]] = \mathbb{E}_{X \sim \mu}[L_{P_X}(s(X))]. \tag{12}$$

We denote by $\mathcal{S}$ the collection of all measurable mappings $s : \mathcal{X} \to \mathfrak{S}_n$, its elements will be referred to as *predictive ranking rules*. As the minimum of the quantity inside the expectation is attained as soon as $s(X)$ is a median for $P_X$, the set of optimal predictive rules can be easily made explicit, as shown by the proposition below.

**Proposition 5** (OPTIMAL ELEMENTS) *The set $\mathcal{S}^*$ of minimizers of the risk (12) is composed of all measurable mappings $s^* : \mathcal{X} \to \mathfrak{S}_n$ such that $s^*(X) \in \mathcal{M}_X$ with probability one, denoting by $\mathcal{M}_x$ the set of median rankings related to distribution $P_x$, $x \in \mathcal{X}$.*

For this reason, the predictive problem formulated above is referred to as *ranking median regression* and its solutions as *conditional median rankings*. It extends the ranking aggregation problem in the sense that $\mathcal{S}^*$ coincides with the set of medians of the marginal distribution $P$ when $\Sigma$ is independent from $X$. Equipped with the notations above, notice incidentally that the minimum prediction error can be written as $\mathcal{R}^* = \mathbb{E}_{X \sim \mu}[L^*_{P_X}]$ and that the risk excess of any $s \in \mathcal{S}$ can be controlled as follows:

$$\mathcal{R}(s) - \mathcal{R}^* \leq \mathbb{E}\left[d\left(s(X),\ s^*(X)\right)\right],$$

for any $s^* \in \mathcal{S}^*$. We assume from now on that $d = d_\tau$. If $P_X \in \mathcal{T}$ with probability one, we almost-surely have $s^*(X) = \sigma^*_{P_X}$ and

$$\mathcal{R}^* = \sum_{i<j} \left\{ 1/2 - \int_{x \in \mathcal{X}} |p_{i,j}(x) - 1/2|\, \mu(dx) \right\},$$

where $p_{i,j}(x) = \mathbb{P}\{\Sigma(i) < \Sigma(j) \mid X = x\}$ for all $i < j$, $x \in \mathcal{X}$. Observe also that in this case, the excess of risk is given by: $\forall s \in \mathcal{S}$,

$$\mathcal{R}(s) - \mathcal{R}^* = \sum_{i<j} \int_{x \in \mathcal{X}} |p_{i,j}(x) - 1/2| \mathbb{I}\{(s(x)(j) - s(x)(i))\,(p_{i,j}(x) - 1/2) < 0\} \mu(dx). \quad (13)$$

The equation above shall play a crucial role in the subsequent fast rate analysis, see Proposition 7's proof in the Appendix.

**Statistical setting.** We assume that we observe $(X_1,\ \Sigma_1)\ \ldots,\ (X_1,\ \Sigma_N)$, $N \geq 1$ i.i.d. copies of the pair $(X,\ \Sigma)$ and, based on these training data, the objective is to build a predictive ranking rule $s$ that nearly minimizes $\mathcal{R}(s)$ over the class $\mathcal{S}$ of measurable mappings $s : \mathcal{X} \to \mathfrak{S}_n$. Of course, the Empirical Risk Minimization (ERM) paradigm encourages to consider solutions of the empirical minimization problem:

$$\min_{s \in \mathcal{S}_0} \widehat{\mathcal{R}}_N(s), \quad (14)$$

where $\mathcal{S}_0$ is a subset of $\mathcal{S}$, supposed to be rich enough for containing approximate versions of elements of $\mathcal{S}^*$ (i.e. so that $\inf_{s \in \mathcal{S}_0} \mathcal{R}(s) - \mathcal{R}^*$ is 'small') and ideally appropriate for continuous or greedy optimization, and

$$\widehat{\mathcal{R}}_N(s) = \frac{1}{N} \sum_{i=1}^N d_\tau(s(X_i),\ \Sigma_i) \quad (15)$$

is a statistical version of (12) based on the $(X_i, \Sigma_i)$'s. Extending those established by Korba et al. (2017) in the context of ranking aggregation, statistical results describing the generalization capacity of minimizers of (15) can be established under classic complexity assumptions for the class $\mathcal{S}_0$, such as the following one (observe incidentally that it is fulfilled by the class of ranking rules output by the algorithm described in subsection 4.3, *cf* Remark 20).

**Assumption 1** *For all $i < j$, the collection of sets*

$$\{\{x \in \mathcal{X} :\ s(x)(i) - s(x)(j) > 0\} :\ s \in \mathcal{S}_0\} \cup \{\{x \in \mathcal{X} :\ s(x)(i) - s(x)(j) < 0\} :\ s \in \mathcal{S}_0\}$$

*is of finite* VC *dimension $V < \infty$.*

**Proposition 6** *Suppose that the class $\mathcal{S}_0$ fulfills Assumption 1. Let $\widehat{s}_N$ be any minimizer of the empirical risk (15) over $\mathcal{S}_0$. For any $\delta \in (0,1)$, we have with probability at least $1-\delta$: $\forall N \geq 1$,*

$$\mathcal{R}(\widehat{s}_N) - \mathcal{R}^* \leq C\sqrt{\frac{V \log(n(n-1)/(2\delta))}{N}} + \left\{\mathcal{R}^* - \inf_{s \in \mathcal{S}_0} \mathcal{R}(s)\right\}, \qquad (16)$$

*where $C < +\infty$ is a universal constant.*

Refer to the Appendix for the technical proof. It is also established there that the rate bound $O_{\mathbb{P}}(1/\sqrt{N})$ is sharp in the minimax sense, see Remark 15.

**Faster learning rates.** As recalled in Section 2, it is proved that rates of convergence for the excess of risk of empirical Kemeny medians can be much faster than $O_{\mathbb{P}}(1/\sqrt{N})$ under transitivity and a certain noise condition (8), see Theorem 3. We now introduce the following hypothesis, involved in the subsequent analysis.

**Assumption 2** *For all $x \in \mathcal{X}$, $P_x \in \mathcal{T}$ and $H = \inf_{x \in \mathcal{X}} \min_{i<j} |p_{i,j}(x) - 1/2| > 0$.*

This condition is checked in many situations, including most conditional parametric models (see Remark 13 in Korba et al. (2017)), and generalizes condition (8), which corresponds to Assumption 2 when $X$ and $\Sigma$ are independent. The result stated below reveals that a similar fast rate phenomenon occurs for minimizers of the empirical risk (15) if Assumption 2 is satisfied. Refer to the Appendix for the technical proof. Since the goal is to give the main ideas, it is assumed for simplicity that the class $\mathcal{S}_0$ is of finite cardinality and that the optimal ranking median regression rule $\sigma^*_{P_x}$ belongs to it.

**Proposition 7** *Suppose that Assumption 2 is fulfilled, that the cardinality of class $\mathcal{S}_0$ is equal to $C < +\infty$ and that the unique true risk minimizer $s^*(x) = \sigma^*_{P_x}$ belongs to $\mathcal{S}_0$. Let $\widehat{s}_N$ be any minimizer of the empirical risk (15) over $\mathcal{S}_0$. For any $\delta \in (0,1)$, we have with probability at least $1-\delta$:*

$$\mathcal{R}(\widehat{s}_N) - \mathcal{R}^* \leq \left(\frac{n(n-1)}{2H}\right) \times \frac{\log(C/\delta)}{N}. \qquad (17)$$

Regarding the minimization problem (14), attention should be paid to the fact that, in contrast to usual (median/quantile) regression, the set $\mathcal{S}$ of predictive ranking rules is not a vector space, which makes the design of practical optimization strategies challenging and the implementation of certain methods, based on (forward stagewise) additive modelling for instance, unfeasible (unless the constraint that predictive rules take their values in $\mathfrak{S}_n$ is relaxed, see Clémençon and Jakubowicz (2010) or Fogel et al. (2013)). If $\mu$ is continuous (the $X_i$'s are pairwise distinct), it is always possible to find $s \in \mathcal{S}$ such that $\widehat{R}_N(s) = 0$ and model selection/regularization issues (*i.e.* choosing an appropriate class $\mathcal{S}_0$) are crucial. In contrast, if $X$ takes discrete values only (corresponding to possible requests in a search engine for instance, like in the usual 'learning to order' setting), in the set $\{1, \ldots, K\}$ with $K \geq 1$ say, the problem (14) boils down to solving *independently* $K$ empirical ranking median problems. However, $K$ may be large and it may be relevant to use some regularization procedure accounting for the possible amount of similarity shared by certain requests/tasks, adding some penalization term to (15). The approach to ranking median regression we develop in this paper, close in spirit to adaptive approximation methods, relies on the concept of local learning and permits to derive practical procedures for building piecewise constant ranking rules (the complexity of the related classes $\mathcal{S}_0$ can be naturally described by the number of constant pieces involved in the predictive rules) from efficient (approximate) Kemeny aggregation (such as that investigated in subsection 2.1), when implemented at a local level. The first method is a version of the popular nearest-neighbor technique, tailored to the ranking median regression setup, while the second algorithm is inspired by the CART algorithm and extends that introduced in Yu et al. (2010), see also Chapter 10 in Alvo and Yu (2014).

## 4. Local Consensus Methods for Ranking Median Regression

We start here with introducing notations to describe the class of piecewise constant ranking rules and explore next approximation of a given ranking rule $s(x)$ by elements of this class, based on a local version of the concept of ranking median recalled in the previous section. Two strategies are next investigated in order to generate adaptively a partition tailored to the training data and yielding a ranking rule with nearly minimum predictive error. Throughout this section, for any measurable set $\mathcal{C} \subset \mathcal{X}$ weighted by $\mu(x)$, the conditional distribution of $\Sigma$ given $X \in \mathcal{C}$ is denoted by $P_\mathcal{C}$. When it belongs to $\mathcal{T}$, the unique median of distribution $P_\mathcal{C}$ is denoted by $\sigma_\mathcal{C}^*$ and referred to as the local median on region $\mathcal{C}$.

### 4.1 Piecewise Constant Predictive Ranking Rules and Local Consensus

Let $\mathcal{P}$ be a partition of $\mathcal{X}$ composed of $K \geq 1$ cells $\mathcal{C}_1$, ..., $\mathcal{C}_K$ (i.e. the $\mathcal{C}_k$'s are pairwise disjoint and their union is the whole feature space $\mathcal{X}$). Suppose in addition that $\mu(\mathcal{C}_k) > 0$ for $k = 1$, ..., $K$. Using the natural embedding $\mathfrak{S}_n \subset \mathbb{R}^n$, any ranking rule $s \in \mathcal{S}$ that is constant on each subset $\mathcal{C}_k$ can be written as

$$s_{\mathcal{P},\bar{\sigma}}(x) = \sum_{k=1}^{K} \sigma_k \cdot \mathbb{I}\{x \in \mathcal{C}_k\}, \tag{18}$$

where $\bar{\sigma} = (\sigma_1, \ldots, \sigma_K)$ is a collection of $K$ permutations. We denote by $\mathcal{S}_\mathcal{P}$ the collection of all ranking rules that are constant on each cell of $\mathcal{P}$. Notice that $\#\mathcal{S}_\mathcal{P} = K \times n!$.

**Local Ranking Medians.** The following result describes the most accurate ranking median regression function in this class. The values it takes correspond to *local Kemeny medians*, i.e. medians of the $P_{\mathcal{C}_k}$'s. The proof is straightforward and postponed to the Appendix.

**Proposition 8** *The set $\mathcal{S}_\mathcal{P}^*$ of solutions of the risk minimization problem $\min_{s \in \mathcal{S}_\mathcal{P}} \mathcal{R}(s)$ is composed of all scoring functions $s_{\mathcal{P},\bar{\sigma}}(x)$ such that, for all $k \in \{1, \ldots, K\}$, the permutation $\sigma_k$ is a Kemeny median of distribution $P_{\mathcal{C}_k}$ and*

$$\min_{s \in \mathcal{S}_\mathcal{P}} \mathcal{R}(s) = \sum_{k=1}^{K} \mu(\mathcal{C}_k) L_{P_{\mathcal{C}_k}}^*.$$

If $P_{\mathcal{C}_k} \in \mathcal{T}$ for $1 \leq k \leq K$, there exists a unique risk minimizer over class $\mathcal{S}_\mathcal{P}$ given by: $\forall x \in \mathcal{X}$,

$$s_\mathcal{P}^*(x) = \sum_{k=1}^K \sigma_{P_{\mathcal{C}_k}}^* \cdot \mathbb{I}\{x \in \mathcal{C}_k\}. \tag{19}$$

Attention should be paid to the fact that the bound

$$\min_{s \in \mathcal{S}_\mathcal{P}} \mathcal{R}(s) - \mathcal{R}^* \leq \inf_{s \in \mathcal{S}_\mathcal{P}} \mathbb{E}_X \left[ d_\tau \left( s^*(X), s(X) \right) \right] \tag{20}$$

valid for all $s^* \in \mathcal{S}^*$, shows in particular that the bias of ERM over the class $\mathcal{S}_\mathcal{P}$ can be controlled by the approximation rate of optimal ranking rules by elements of $\mathcal{S}_\mathcal{P}$ when error is measured by the integrated Kendall $\tau$ distance and $X$'s marginal distribution, $\mu(x)$ namely, is the integration measure.

**Approximation.** We now investigate to what extent ranking median regression functions $s^*(x)$ can be well approximated by predictive rules of the form (18). We assume that $\mathcal{X} \subset \mathbb{R}^d$ with $d \geq 1$ and denote by $||.||$ any norm on $\mathbb{R}^d$. The following hypothesis is a classic smoothness assumption on the conditional pairwise probabilities.

**Assumption 3** *For all $1 \leq i < j \leq n$, the mapping $x \in \mathcal{X} \mapsto p_{i,j}(x)$ is Lipschitz, i.e. there exists $M < \infty$ such that:*

$$\forall (x, x') \in \mathcal{X}^2, \ \sum_{i<j} |p_{i,j}(x) - p_{i,j}(x')| \leq M \cdot ||x - x'||. \tag{21}$$

The following result shows that, under the assumptions above, the optimal prediction rule $\sigma_{P_X}^*$ can be accurately approximated by (19), provided that the regions $\mathcal{C}_k$ are 'small' enough.

**Theorem 9** *Suppose that $P_x \in \mathcal{T}$ for all $x \in \mathcal{X}$ and that Assumption 3 is fulfilled. Then, we have: $\forall s_\mathcal{P} \in \mathcal{S}_\mathcal{P}^*$.*

$$\mathcal{R}(s_\mathcal{P}) - \mathcal{R}^* \leq M \cdot \delta_\mathcal{P}, \tag{22}$$

*where $\delta_\mathcal{P} = \max_{\mathcal{C} \in \mathcal{P}} \sup_{(x,x') \in \mathcal{C}^2} ||x - x'||$ is the maximal diameter of $\mathcal{P}$'s cells. Hence, if $(\mathcal{P}_m)_{m \geq 1}$ is a sequence of partitions of $\mathcal{X}$ such that $\delta_{\mathcal{P}_m} \to 0$ as $m$ tends to infinity, then $\mathcal{R}(s_{\mathcal{P}_m}) \to \mathcal{R}^*$ as $m \to \infty$.*

*Suppose in addition that Assumption 2 is fulfilled and that $P_\mathcal{C} \in \mathcal{T}$ for all $\mathcal{C} \in \mathcal{P}$. Then, we have:*

$$\mathbb{E} \left[ d_\tau \left( \sigma_{P_X}^*, s_\mathcal{P}^*(X) \right) \right] \leq \sup_{x \in \mathcal{X}} d_\tau \left( \sigma_{P_x}^*, s_\mathcal{P}^*(x) \right) \leq (M/H) \cdot \delta_\mathcal{P}. \tag{23}$$

The upper bounds above reflect the fact that the the smaller the Lipschitz constant $M$, the easier the ranking median regression problem and that the larger the quantity $H$, the easier the recovery of the optimal RMR rule. In the Appendix, examples of distributions $(\mu(dx), P_x)$ satisfying Assumptions 2-3 both at the same time are given.

**Remark 10** (ON LEARNING RATES) *For simplicity, assume that $\mathcal{X} = [0,1]^d$ and that $\mathcal{P}_m$ is a partition with $m^d$ cells with diameter less than $C \times 1/m$ each, where $C$ is a constant. Provided the assumptions it stipulates are fulfilled, Theorem 9 shows that the bias of the ERM method over the class $\mathcal{S}_{\mathcal{P}_m}$ is of order $1/m$. Combined with Proposition 6, choosing $m \sim \sqrt{N}$ gives a nearly optimal learning rate, of order $O_\mathbb{P}((\log N)/N)$ namely.*

**Remark 11** (ON SMOOTHNESS ASSUMPTIONS) *We point out that the analysis above could be naturally refined, insofar as the accuracy of a piecewise constant median ranking regression rule is actually controlled by its capacity to approximate an optimal rule $s^*(x)$ in the $\mu$-integrated Kendall $\tau$ sense, as shown by Eq. (20). Like in Binev et al. (2005) for distribution-free regression, learning rates for ranking median regression could be investigated under the assumption that $s^*$ belongs to a certain smoothness class defined in terms of approximation rate, specifying the decay rate of $\inf_{s \in \mathcal{S}_m} \mathbb{E}[d_\tau(s^*(X), s(X))]$ for a certain sequence $(\mathcal{S}_m)_{m \geq 1}$ of classes of piecewise constant ranking rules. This is beyond the scope of the present paper and will be the subject of future work.*

The next result, proved in the Appendix section, states a very general consistency theorem for a wide class of RMR rules based on data-based partitioning, in the spirit of Lugosi and Nobel (1996) for classification. For simplicity's sake, we assume that $\mathcal{X}$ is compact, equal to $[0,1]^d$ say. Let $N \geq 1$, a $N$-sample partitioning rule $\pi_N$ maps any possible training sample $\mathcal{D}_N = ((x_1, \sigma_1), \ldots, (x_N, \sigma_N)) \in (\mathcal{X} \times \mathfrak{S}_n)^N$ to a partition $\pi_N(\mathcal{D}_N)$ of $[0,1]^d$ composed of borelian cells. The associated collection of partitions is denoted by $\mathcal{F}_N = \{\pi_N(\mathcal{D}_N) : \mathcal{D}_N \in (\mathcal{X} \times \mathfrak{S}_n)^N\}$. As in Lugosi and Nobel (1996), the complexity of $\mathcal{F}_N$ is measured by the $N$-order shatter coefficient of the class of sets that can be obtained as unions of cells of a partition in $\mathcal{F}_N$, denoted by $\Delta_N(\mathcal{F}_N)$. An estimate of this quantity can be found in *e.g.* Chapter 21 of Devroye et al. (1996) for various data-dependent partitioning rules (including the recursive partitioning scheme described in subsection 4.3, when implemented with axis-parallel splits). When $\pi_N$ is applied to a training sample $\mathcal{D}_N$, it produces a partition $\mathcal{P}_N = \pi_N(\mathcal{D}_N)$ (that is random in the sense that it depends on $\mathcal{D}_N$) associated with a RMR prediction rule: $\forall x \in \mathcal{X}$,

$$s_N(x) = \sum_{\mathcal{C} \in \mathcal{P}_N} \widehat{\sigma}_\mathcal{C} \cdot \mathbb{I}\{x \in \mathcal{C}\} \tag{24}$$

where $\widehat{\sigma}_\mathcal{C}$ denotes a Kemeny median of the empirical version of $\Sigma$'s distribution given $X \in \mathcal{C}$, $\widehat{P}_\mathcal{C} = (1/N_\mathcal{C}) \sum_{i: X_i \in \mathcal{C}} \delta_{\Sigma_i}$ with $N_\mathcal{C} = \sum_i \mathbb{I}\{X_i \in \mathcal{C}\}$ and the convention $0/0 = 0$, for any measurable set $\mathcal{C}$ s.t. $\mu(\mathcal{C}) > 0$. Notice that, although we have $\widehat{\sigma}_\mathcal{C} = \sigma^*_{\widehat{P}_\mathcal{C}}$ if $\widehat{P}_\mathcal{C} \in \mathcal{T}$, the rule 24 is somehow theoretical, since the way the Kemeny medians $\widehat{\sigma}_\mathcal{C}$ are obtained is not specified in general. Alternatively, using the notations of subsection 2.1, one may consider the RMR rule

$$\widetilde{s}_N(x) = \sum_{\mathcal{C} \in \mathcal{P}_N} \widetilde{\sigma}^*_{\widehat{P}_\mathcal{C}} \cdot \mathbb{I}\{x \in \mathcal{C}\}, \tag{25}$$

which takes values that are not necessarily local empirical Kemeny medians but can always be easily computed. Observe incidentally that, for any $\mathcal{C} \in \mathcal{P}_N$ s.t. $\widehat{P}_\mathcal{C} \in \mathcal{T}$, we have $\widetilde{s}_N(x) = s_N(x)$ for all $x \in \mathcal{C}$. The theorem below establishes the consistency of these RMR rules in situations where the diameter of the cells of the data-dependent partition and their $\mu$-measure decay to zero but not too fast, with respect to the rate at which the quantity $\sqrt{N/\log(\Delta_n(\mathcal{F}_N))}$ increases.

**Theorem 12** *Let $(\pi_1, \pi_2, \ldots)$ be a fixed sequence of partitioning rules and for each $N$ let $\mathcal{F}_N$ be the collection of partitions associated with the $N$-sample partitioning rule $\pi_N$. Suppose that $P_x \in \mathcal{T}$ for all $x \in \mathcal{X}$ and that Assumption 3 is satisfied. Assume also that the conditions below are fulfilled:*

(i) $\lim_{n\to\infty} \log(\Delta_N(\mathcal{F}_N))/N = 0$,

(ii) *we have* $\delta_{\mathcal{P}_N} \to 0$ *in probability as* $N \to \infty$ *and*

$$1/\kappa_N = o_{\mathbb{P}}(\sqrt{N/\log \Delta_N(\mathcal{F}_N)}) \text{ as } N \to \infty,$$

*where* $\kappa_N = \inf\{\mu(\mathcal{C}) : \mathcal{C} \in \mathcal{P}_N\}$.

*Then any RMR rule* $s_N$ *of the form* (24) *is consistent, i.e.* $\mathcal{R}(s_N) \to \mathcal{R}^*$ *in probability as* $N \to \infty$.
*Suppose in addition that Assumption 2 is satisfied. Then, the RMR rule* $\tilde{s}_N(x)$ *given by* (25) *is also consistent.*

The next section presents two approaches for building a partition $\mathcal{P}$ of the predictor variable space in a data-driven fashion. The first method is a version of the nearest neighbor methods tailored to ranking median regression, whereas the second algorithm constructs $\mathcal{P}$ recursively, depending on the local variability of the $\Sigma_i$'s, and scales with the dimension of the input space.

### 4.2 Nearest-Neighbor Rules for Ranking Median Regression

Fix $k \in \{1, \ldots, N\}$ and a query point $x \in \mathcal{X}$. The $k$-nearest neighbor RMR rule prediction $s_{k,N}(x)$ is obtained as follows. Sort the training data $(X_1, \Sigma_1), \ldots, (X_n, \Sigma_n)$ by increasing order of the distance to $x$, measured, for simplicity, by $\|X_i - x\|$ for a certain norm chosen on $\mathcal{X} \subset \mathbb{R}^d$ say: $\|X_{(1,N)} - x\| \leq \ldots \leq \|X_{(N,N)} - x\|$. Consider next the empirical distribution calculated using the $k$ training points closest to $x$

$$\widehat{P}(x) = \frac{1}{k} \sum_{l=1}^{k} \delta_{\Sigma_{(l,N)}} \tag{26}$$

and then set

$$s_{k,N}(x) = \sigma_{\widehat{P}(x)}, \tag{27}$$

where $\sigma_{\widehat{P}(x)}$ is a Kemeny median of distribution (26). Alternatively, one may compute next the pseudo-empirical Kemeny median, as described in subsection 2.1, yielding the $k$-NN prediction at $x$:

$$\tilde{s}_{k,N}(x) = \tilde{\sigma}^*_{\widehat{P}(x)}. \tag{28}$$

Observe incidentally that $s_{k,N}(x) = \tilde{s}_{k,N}(x)$ when $\widehat{P}(x)$ is strictly stochastically transitive. The result stated below provides an upper bound for the expected risk excess of the RMR rules (27) and (28), which reflects the usual bias/variance trade-off ruled by $k$ for fixed $N$ and asymptotically vanishes as soon as $k \to \infty$ as $N \to \infty$ such that $k = o(N)$. Notice incidentally that the choice $k \sim N^{2/(d+2)}$ yields the asymptotically optimal upper bound, of order $N^{-1/(2+d)}$.

**Theorem 13** *Suppose that Assumption 3 is fulfilled, that the r.v. $X$ is bounded and $d \geq 3$. Then, we have: $\forall N \geq 1$, $\forall k \in \{1, \ldots, N\}$,*

$$\mathbb{E}\left[\mathcal{R}(s_{k,N}) - \mathcal{R}^*\right] \leq \frac{n(n-1)}{2} \left( \frac{1}{\sqrt{k}} + 2\sqrt{c_1} M \left( \frac{k}{N} \right)^{1/d} \right) \tag{29}$$

where $c_1$ is a constant which only depends on $\mu$'s support.

Suppose in addition that Assumption 2 is satisfied. We then have: $\forall N \geq 1$, $\forall k \in \{1, \ldots, N\}$,

$$\mathbb{E}\left[\mathcal{R}(\widetilde{s}_{k,N}) - \mathcal{R}^*\right] \leq \frac{n(n-1)}{2} \left( \frac{1}{\sqrt{k}} + 2\sqrt{c_1} M \left( \frac{k}{N} \right)^{1/d} \right) (1 + n(n-1)/(4H)). \tag{30}$$

Refer to the Appendix for the technical proof. In addition, for $d \leq 2$ the rate stated in Theorem 13 still holds true, under additional conditions on $\mu$, see the Appendix for further details. In practice, as for nearest-neighbor methods in classification/regression, the success of the technique above for fixed $N$ highly depends on the number $k$ of neighbors involved in the computation of the local prediction. The latter can be picked by means of classic model selection methods, based on data segmentation/resampling techniques. It may also crucially depend on the distance chosen (which could be learned from the data as well, see *e.g.* Bellet et al. (2013)) and/or appropriate preprocessing stages, see *e.g.* the discussion in chapter 13 of Friedman et al. (2002)). The implementation of this simple local method for ranking median regression does not require to explicit the underlying partition but is classically confronted with the curse of dimensionality. The next subsection explains how another local method, based on the popular tree induction heuristic, scales with the dimension of the input space by contrast. Due to space limitations, extensions of other data-dependent partitioning methods, such as those investigated in Chapter 21 of Devroye et al. (1996) for instance, to local RMR are left to the reader.

### 4.3 Recursive Partitioning - The CRIT algorithm

We now describe an iterative scheme for building an appropriate tree-structured partition $\mathcal{P}$, adaptively from the training data. Whereas the splitting criterion in most recursive partitioning methods is heuristically motivated (see Friedman (1997)), the local learning method we describe below relies on the Empirical Risk Minimization principle formulated in Section 3, so as to build by refinement a partition $\mathcal{P}$ based on a training sample $\mathcal{D}_N = \{(\Sigma_1, \mathcal{X}_1), \ldots, (\Sigma_N, X_N)\}$ so that, on each cell $\mathcal{C}$ of $\mathcal{P}$, the $\Sigma_i$'s lying in it exhibit a small variability in the Kendall $\tau$ sense and, consequently, may be accurately approximated by a local Kemeny median. As shown below, the local variability measure we consider can be connected to the local ranking median regression risk (see Eq. (34)) and leads to exactly the same node impurity measure as in the tree induction method proposed in Yu et al. (2010), see Remark 14. The algorithm described below differs from it in the method we use to compute the local predictions. More precisely, the goal pursued is to construct recursively a piecewise constant ranking rule associated to a partition $\mathcal{P}$, $s_{\mathcal{P}}(x) = \sum_{\mathcal{C} \in \mathcal{P}} \sigma_{\mathcal{C}} \cdot \mathbb{I}\{x \in \mathcal{C}\}$, with minimum empirical risk

$$\widehat{L}_N(s_\mathcal{P}) = \sum_{\mathcal{C} \in \mathcal{P}} \widehat{\mu}_N(\mathcal{C}) L_{\widehat{P}_\mathcal{C}}(\sigma_\mathcal{C}), \tag{31}$$

where $\widehat{\mu}_N = (1/N) \sum_{k=1}^N \delta_{X_k}$ is the empirical measure of the $X_k$'s. The partition $\mathcal{P}$ being fixed, as noticed in Proposition 8, the quantity (31) is minimum when $\sigma_\mathcal{C}$ is a Kemeny median of $\widehat{P}_\mathcal{C}$ for all $\mathcal{C} \in \mathcal{P}$. It is then equal to

$$\min_{s \in \mathcal{S}_\mathcal{P}} \widehat{L}_N(s) = \sum_{\mathcal{C} \in \mathcal{P}} \widehat{\mu}_N(\mathcal{C}) L^*_{\widehat{P}_\mathcal{C}}. \tag{32}$$

Except in the case where the intra-cell empirical distributions $\widehat{P}_\mathcal{C}$'s are all stochastically transitive (each $L^*_{\widehat{P}_\mathcal{C}}$ can be then computed using formula (6)), computing (32) at each recursion of the algorithm can be very expensive, since it involves the computation of a Kemeny median within each cell $\mathcal{C}$. We propose to measure instead the accuracy of the current partition by the quantity

$$\widehat{\gamma}_\mathcal{P} = \sum_{\mathcal{C} \in \mathcal{P}} \widehat{\mu}_N(\mathcal{C}) \gamma_{\widehat{P}_\mathcal{C}}, \tag{33}$$

which satisfies the double inequality (see Remark 2)

$$\widehat{\gamma}_\mathcal{P} \leq \min_{s \in \mathcal{S}_\mathcal{P}} \widehat{L}_N(s) \leq 2\widehat{\gamma}_\mathcal{P}, \tag{34}$$

and whose computation is straightforward: $\forall \mathcal{C} \in \mathcal{P}$,

$$\gamma_{\widehat{P}_\mathcal{C}} = \frac{1}{2} \sum_{i<j} \widehat{p}_{i,j}(\mathcal{C}) (1 - \widehat{p}_{i,j}(\mathcal{C})), \tag{35}$$

where $\widehat{p}_{i,j}(\mathcal{C}) = (1/N_\mathcal{C}) \sum_{k: X_k \in \mathcal{C}} \mathbb{I}\{\Sigma_k(i) < \Sigma_k(j)\}$, $i < j$, denote the local pairwise empirical probabilities. A ranking median regression tree of maximal depth $J \geq 0$ is grown as follows. One starts from the root node $\mathcal{C}_{0,0} = \mathcal{X}$. At depth level $0 \leq j < J$, any cell $\mathcal{C}_{j,k}$, $0 \leq k < 2^j$ shall be split into two (disjoint) subsets $\mathcal{C}_{j+1,2k}$ and $\mathcal{C}_{j+1,2k+1}$, respectively identified as the left and right children of the interior leaf $(j,k)$ of the ranking median regression tree, according to the following *splitting rule*.

**Splitting rule.** For any candidate left child $\mathcal{C} \subset \mathcal{C}_{j,k}$, picked in a class $\mathcal{G}$ of 'admissible' subsets (see Remark 20), the relevance of the split $\mathcal{C}_{j,k} = \mathcal{C} \cup (\mathcal{C}_{j,k} \setminus \mathcal{C})$ is naturally evaluated through the quantity:

$$\Lambda_{j,k}(\mathcal{C}) \stackrel{def}{=} \widehat{\mu}_N(\mathcal{C}) \gamma_{\widehat{P}_\mathcal{C}} + \widehat{\mu}_N(\mathcal{C}_{j,k} \setminus \mathcal{C}) \gamma_{\widehat{P}_{\mathcal{C}_{j,k} \setminus \mathcal{C}}}. \tag{36}$$

The determination of the splitting thus consists in computing a solution $\mathcal{C}_{j+1,2k}$ of the optimization problem

$$\min_{\mathcal{C} \in \mathcal{G}, \mathcal{C} \subset \mathcal{C}_{j,k}} \Lambda_{j,k}(\mathcal{C}) \tag{37}$$

As explained in the Appendix, an appropriate choice for class $\mathcal{G}$ permits to solve exactly the optimization problem very efficiently, in a greedy fashion.

**Local medians.** The consensus ranking regression tree is grown until depth $J$ and on each terminal leave $\mathcal{C}_{J,l}$, $0 \leq l < 2^J$, one computes the local Kemeny median estimate by means of the best strictly stochastically transitive approximation method investigated in subsection 2.1

$$\sigma_{J,l}^* \stackrel{def}{=} \widetilde{\sigma}_{\widehat{P}_{\mathcal{C}_{J,l}}}^*. \tag{38}$$

If $\widehat{P}_{\mathcal{C}_{J,l}} \in \mathcal{T}$, $\sigma_{J,l}^*$ is straightforwardly obtained from formula (7) and otherwise, one uses the pseudo-empirical Kemeny median described in subsection 2.1. The ranking median regression rule related to the binary tree $T_{2^J}$ thus constructed is given by:

$$s_{T_{2^J}}^*(x) = \sum_{l=0}^{2^J - 1} \sigma_{J,l}^* \mathbb{I}\{x \in \mathcal{C}_{J,l}\}. \tag{39}$$

Its training prediction error is equal to $\widehat{L}_N(s_{T_{2^J}}^*)$, while the training accuracy measure of the final partition is given by

$$\widehat{\gamma}_{T_{2^J}} = \sum_{l=0}^{2^J - 1} \widehat{\mu}_N(\mathcal{C}_{J,l}) \gamma_{\widehat{P}_{\mathcal{C}_{J,l}}}. \tag{40}$$

**Remark 14** *We point out that the impurity measure (33) corresponds (up to a constant factor) to that considered in Yu et al. (2010), where it is referred to as the pairwise Gini criterion. Borrowing their notation, one may indeed write: for any measurable $\mathcal{C} \subset \mathcal{X}$, $i_w^{(2)}(\mathcal{C}) = 8/(n(n-1)) \times \gamma_{\widehat{P}_{\mathcal{C}}}$.*

The tree growing stage is summarized in the Appendix, as well as the pruning procedure generally following it to avoid overfitting. Additional comments on the advantages of this method regarding interpretability and computational feasibility can also be found in the Appendix, together with a preliminary analysis of a specific bootstrap aggregation technique that can remedy the instability of such a hierarchical predictive method.

## 5. Numerical Experiments

For illustration purpose, experimental results based on simulated/real data are displayed.

### 5.1 Results on simulated data

Here, datasets of full rankings on $n$ items are generated according to two explanatory variables. We carried out several experiments by varying the number of items ($n = 3, 5, 8$) and the nature of the features. In Setting 1, both features are numerical; in Setting 2, one is numerical and the other categorical, while, in Setting 3, both are categorical. For a fixed setting, a partition $\mathcal{P}$ of $\mathcal{X}$ composed of $K$ cells $\mathcal{C}_1, \ldots, \mathcal{C}_K$ is fixed. In each trial, $K$ permutations $\sigma_1, \ldots, \sigma_K$ (which can be arbitrarily close) are generated, as well as three datasets of $N$ samples, where on each cell $\mathcal{C}_k$: the first one is constant (all samples are equal to $\sigma_k$), and the two others are noisy versions of the first one, where the samples follow a Mallows distribution (see Mallows (1957)) centered on $\sigma_k$ with dispersion parameter $\phi$. We recall that the greater the dispersion parameter $\phi$, the spikiest the distribution (and closest to piecewise constant). We choose $K=6$ and $N=1000$. In each trial, the dataset is divided into a training set (70%) and a test set (30%). Concerning the CRIT algorithm, since the true partition is known and is of depth 3, the maximum depth is set to 3 and the minimum size in a leaf is set to the number of samples in the training set divided by 10. For the k-NN algorithm, the number of neighbors $k$ is fixed to 5. The baseline model to which we compare our algorithms is the following: on the train set, we fit a K-means (with $K=6$), train a Plackett-Luce model on each cluster and assign the mode of this learnt distribution as the center ranking of the cluster. For each configuration (number of items, characteristics of feature and distribution of the dataset), the empirical risk (see 3, denoted as $\widehat{\mathcal{R}}_N(s)$) is averaged on 50 repetitions of the experiment. Results of the k-NN algorithm (indicated with a star *), of the CRIT algorithm (indicated with two stars **) and of the baseline model (between parenthesis) on the various configurations are provided in Table 1. They show that the methods we develop recover the true partition of the data, insofar as the underlying distribution can be well approximated by a piecewise constant function ($\phi \geq 2$ for instance in our simulations).

### 5.2 Analysis of GSS data on job value preferences

We test our algorithm on the full rankings dataset which was obtained by the US General Social Survey (GSS) and which is already used in Alvo and Yu (2014). This multidimensional survey collects across years socio-demographic attributes and answers of respondents to numerous questions, including societal opinions. In particular, participants were asked to rank in order of preference five aspects about a job: "high income", "no danger of being fired", "short working hours", "chances for advancement", and "work important and gives a feeling of accomplishment". The dataset we consider contains 18544 samples collected between 1973 and 2014. As in Alvo and Yu (2014), for each individual, we consider eight individual attributes (sex, race, birth cohort, highest educational degree attained, family income, marital status, number of children that the respondent ever had, and household size) and three properties of work conditions (working status, employment status, and occupation). We repeat 10 times the experiment, each time by randomly dividing the dataset in a training set (70%) and a test set (30%), and training the model. The results are stable among the experiments. Concerning the k-NN algorithm, we obtain an average empirical risk of 3.235 (for $k = 3$). For the CRIT algorithm, we obtain an average empirical risk of 2.763 (recall that the maximum Kendall distance is 10) and splits coherent with the analysis in Alvo and Yu (2014): the first splitting variable is occupation (managerial, professional, sales workers and related vs services, natural resources, production, construction and transportation occupations), then at the second level the race is the most important factor in both groups (black respondents vs others in the first group, white respondents vs others in the second group). At the lower level the degree obtained seems to play an important role (higher than high school, or higher than bachelor's for example in some groups); then other discriminating variables among lower levels are birth cohort, family income or working status.

## 6. Conclusion and Perspectives

The contribution of this article is twofold. The problem of learning to predict preferences, expressed in the form of a permutation, in a supervised setting is formulated and investigated in a rigorous probabilistic framework (optimal elements, learning rate bounds, bias analysis), extending that recently developed for statistical Kemeny ranking aggregation in Korba et al. (2017). Based on this formulation, it is also shown that predictive methods based on the concept of local Kemeny consensus, variants of nearest-neighbor and tree-induction methods namely, are well-suited for this learning task. This is justified by approximation theoretic arguments and algorithmic simplicity/efficiency both at the same time and illustrated by numerical experiments. Whereas the ranking median regression problem is motivated by many applications in our era of recommender systems and personalized customer services, the output variable may take the form of an *incomplete ranking* rather than a full ranking in many situations. Extension of the results to this more general framework, following in the footsteps of the multiresolution analysis approach developed in Sibony et al. (2015), will be the subject of further research.

## Appendix

**Proof of Lemma 4**

Observe first that: $\forall \sigma \in \mathfrak{S}_n$,

$$L_{P'}(\sigma) = \sum_{i<j} p'_{i,j}\mathbb{I}\{\sigma(i) > \sigma(j)\} + \sum_{i<j}(1 - p'_{i,j})\mathbb{I}\{\sigma(i) < \sigma(j)\}. \tag{41}$$

We deduce from the equality above, applied twice, that: $\forall \sigma \in \mathfrak{S}_n$,

$$|L_{P'}(\sigma)) - L_{P''}(\sigma)| \leq \sum_{i<j} \left| p'_{i,j} - p''_{i,j} \right|. \tag{42}$$

Hence, we may write:

$$L^*_{P'} = \inf_{\sigma \in \mathfrak{S}_n} L_{P'}(\sigma) \leq L_{P'}(\sigma_{P''}) = L_{P''}(\sigma_{P''}) + (L_{P'}(\sigma_{P''}) - L_{P''}(\sigma_{P''}))$$

$$\leq L^*_{P''} + \sum_{i<j} \left| p'_{i,j} - p''_{i,j} \right|.$$

In a similar fashion, we have $L^*_{P''} \leq L^*_{P'} + \sum_{i<j}|p'_{i,j} - p''_{i,j}|$, which yields assertion (i) when combined with the inequality above.

We turn to (ii) and assume now that both $P'$ and $P''$ belong to $\mathcal{T}$. Let $i < j$. Suppose that $\sigma^*_{P'}(i) < \sigma^*_{P'}(j)$ and $\sigma^*_{P''}(i) > \sigma^*_{P''}(j)$. In this case, we have $p'_{i,j} > 1/2$ and $p''_{i,j} < 1/2$, so that

$$|p'_{i,j} - p''_{i,j}|/h = \left(p'_{i,j} - 1/2\right)/h + \left(1/2 - p''_{i,j}\right)/h \geq 1.$$

More generally, we have

$$\mathbb{I}\left\{\left(\sigma^*_{P'}(i) - \sigma^*_{P'}(j)\right)\left(\sigma^*_{P'}(i) - \sigma^*_{P'}(j)\right) < 0\right\} \leq |p'_{i,j} - p''_{i,j}|/h$$

for all $i < j$. Summing over the pairs $(i,j)$ establishes assertion (ii).

**Proof of Proposition 6**

First, observe that, using the definition of empirical risk minimizers and the union bound, we have with probability one: $\forall N \geq 1$,

$$\mathcal{R}(\widehat{s}_N) - \mathcal{R}^* \leq 2 \sup_{s \in \mathcal{S}_0} \left|\widehat{\mathcal{R}}_N(s) - \mathcal{R}(s)\right| + \left\{\inf_{s \in \mathcal{S}_0} \mathcal{R}(s) - \mathcal{R}^*\right\}$$

$$\leq 2 \sum_{i<j} \sup_{s \in \mathcal{S}_0} \left|\frac{1}{N} \sum_{k=1}^{N} \mathbb{I}\{(\Sigma_k(i) - \Sigma_k(j))(s(X_k)(i) - s(X_k)(j)) < 0\} - \mathcal{R}_{i,j}(s)\right|$$

$$+ \left\{\inf_{s \in \mathcal{S}_0} \mathcal{R}(s) - \mathcal{R}^*\right\},$$

where $\mathcal{R}_{i,j}(s) = \mathbb{P}\{(\Sigma(i) - \Sigma(j))(s(X)(i) - s(X)(j)) < 0\}$ for $i < j$ and $s \in \mathcal{S}$. Since Assumption 1 is satisfied, by virtue of Vapnik-Chervonenkis inequality (see *e.g.* Devroye et al. (1996)), for all $i < j$ and any $\delta \in (0, 1)$, we have with probability at least $1 - \delta$:

$$\sup_{s \in \mathcal{S}_0} \left|\frac{1}{N} \sum_{k=1}^{N} \mathbb{I}\{(\Sigma_k(i) - \Sigma_k(j))(s(X_k)(i) - s(X_k)(j)) < 0\} - \mathcal{R}_{i,j}(s)\right| \leq c\sqrt{V \log(1/\delta)/N}, \tag{43}$$

where $c < +\infty$ is a universal constant. The desired bound then results from the combination of the bound above and the union bound.

**Remark 15** (ON MINIMAXITY) *Observing that, when $X$ and $\Sigma$ are independent, the best predictions are $P$'s Kemeny medians, it follows from Proposition 11 in Korba et al. (2017) that the minimax risk can be bounded by below as follows:*

$$\inf_{s_N} \sup_{\mathcal{L}} \mathbb{E}_{\mathcal{L}}\left[\mathcal{R}_{\mathcal{L}}(s_N) - \mathcal{R}^*_{\mathcal{L}}\right] \geq \frac{1}{16e\sqrt{N}},$$

*where the supremum is taken over all possible probability distributions $\mathcal{L} = \mu(dx) \otimes P_x(d\sigma)$ for $(X, \Sigma)$ (including the independent case) and the minimum is taken over all mappings that map a dataset $(X_1, \Sigma_1)$, ..., $(X_N, \Sigma_N)$ made of independent copies of $(X, \Sigma)$ to a ranking rule in $\mathcal{S}$.*

**Proof of Proposition 7**

The subsequent fast rate analysis mainly relies on the lemma below.

**Lemma 16** *Suppose that Assumption 2 is fulfilled. Let $s \in \mathcal{S}$ and set*

$$Z(s) = \sum_{i<j} \{\mathbb{I}\{(\Sigma(i) - \Sigma(j))(s(X)(i) - s(X)(j)) < 0\} -$$

$$\sum_{i<j} \mathbb{I}\{(\Sigma(i) - \Sigma(j))\left(\sigma^*_{P_X}(i) - \sigma^*_{P_X}(j)\right) < 0\}\}.$$

*Then, we have:*

$$Var(Z(s)) \leq \left(\frac{n(n-1)}{2H}\right) \times (\mathcal{R}(s) - \mathcal{R}^*).$$

**Proof** Recall first that it follows from (5) that, for all $i < j$,
$$\left(\sigma^*_{P_X}(j) - \sigma^*_{P_X}(i)\right)(p_{i,j}(X) - 1/2) > 0.$$

Hence, we have:
$$Var(Z(s)) \leq \frac{n(n-1)}{2} \times \sum_{i<j} Var\left(\mathbb{I}\{(p_{i,j}(X) - 1/2)(s(X)(j) - s(X)(i)) < 0\}\right)$$
$$\leq \frac{n(n-1)}{2} \times \sum_{i<j} \mathbb{E}\left[\mathbb{I}\{(p_{i,j}(X) - 1/2)(s(X)(j) - s(X)(i)) < 0\}\right].$$

In addition, it follows from formula (13) for the risk excess that:
$$\mathcal{R}(s) - \mathcal{R}^* \geq H \times \sum_{i<j} \mathbb{E}\left[\mathbb{I}\{(p_{i,j}(X) - 1/2)(s(X)(j) - s(X)(i)) < 0\}\right].$$

Combined with the previous inequality, this establishes the lemma. ∎

Since the goal is to give the main ideas, we assume for simplicity that the $\mathcal{S}_0$ is of finite cardinality and that the optimal ranking median regression rule $\sigma^*_{P_x}$ belongs to it. Applying Bernstein's inequality to the i.i.d. average $(1/N)\sum_{k=1}^N Z_k(s)$, where
$$Z_k(s) = \sum_{i<j}\{\mathbb{I}\{(\Sigma_k(i) - \Sigma_k(j))(s(X_k)(i) - s(X_k)(j)) < 0\} -$$
$$\sum_{i<j} \mathbb{I}\left\{(\Sigma_k(i) - \Sigma_k(j))\left(\sigma^*_{P_{X_k}}(i) - \sigma^*_{P_{X_k}}(j)\right) < 0\right\}\},$$

for $1 \leq k \leq N$ and the union bound over the ranking rules $s$ in $\mathcal{S}_0$, we obtain that, for all $\delta \in (0,1)$, we have with probability larger than $1 - \delta$: $\forall s \in \mathcal{S}_0$,
$$\mathbb{E}[Z(s)] = \mathcal{R}(s) - \mathcal{R}^* \leq \widehat{\mathcal{R}}_N(s) - \widehat{\mathcal{R}}_N(s^*) + \sqrt{\frac{2Var(Z(s))\log(C/\delta)}{N}} + \frac{4\log(C/\delta)}{3N}.$$

Since $\widehat{\mathcal{R}}_N(\widehat{s}_N) - \widehat{\mathcal{R}}_N(s^*) \leq 0$ by assumption and using the variance control provided by Lemma 16 above, we obtain that, with probability at least $1 - \delta$, we have:
$$\mathcal{R}(\widehat{s}_N) - \mathcal{R}^* \leq \sqrt{\frac{\frac{n(n-1)}{H}(\mathcal{R}(\widehat{s}_N) - \mathcal{R}^*)/H \times \log(C/\delta)}{N}} + \frac{4\log(C/\delta)}{3N}.$$

Finally, solving this inequality in $\mathcal{R}(\widehat{s}_N) - \mathcal{R}^*$ yields the desired result.

**Proof of Proposition 8**

Let $s(x) = \sum_{k=1}^K \sigma_k \mathbb{I}\{x \in \mathcal{C}_k\}$ in $\mathcal{S}_\mathcal{P}$. It suffices to observe that we have
$$\mathcal{R}(s) = \sum_{k=1}^K \mathbb{E}\left[\mathbb{I}\{X \in \mathcal{C}_k\} d_\tau(\sigma_k, \Sigma)\right] = \sum_{k=1}^K \mu(\mathcal{C}_k) \mathbb{E}\left[d_\tau(\sigma_k, \Sigma) \mid X \in \mathcal{C}_k\right], \quad (44)$$

and that each term involved in the summation above is minimum for $\sigma_k \in \mathcal{M}_{P_{\mathcal{C}_k}}$, $1 \leq k \leq K$.

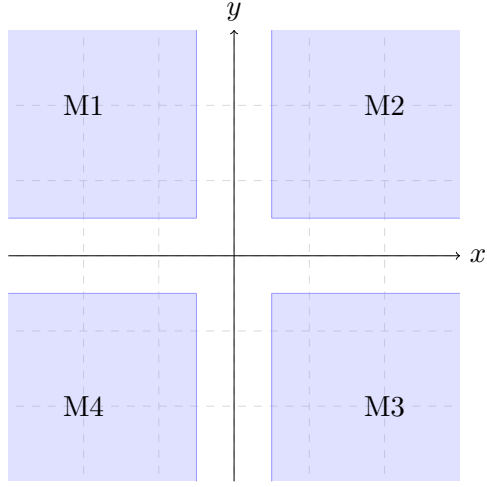

Figure 1: Example of a distribution satisfying Assumptions 2-3 in $\mathbb{R}^2$.

**An example of a distribution satisfying Assumptions 2-3**

Let $d=2$ and $\mathcal{P}$ a partition of $\mathbb{R}^2$ given Figure 1. Suppose that for $x \in \mathcal{X}$, $\mu(x)$ is null outside the filled areas $(M_i)_{i=1,\ldots,4}$, and that on each filled area $M_i$ for $i=1,\ldots,4$, the conditional distribution of $\Sigma$ given $X \in M_i$, denoted by $P_{M_i}$, is a Mallows distribution with parameters $(\pi_i, \phi_i)$. In this case, Assumption 3 is satisfied. Then, if for each $i = 1, \ldots, 4$, $\phi_i \leq (1-2H)/(1+2H)$, we also have Assumption 2 verified.

**Proof of Theorem 9**

Consider $s_\mathcal{P}(x) = \sum_{\mathcal{C} \in \mathcal{P}} \sigma_{P_\mathcal{C}} \cdot \mathbb{I}\{x \in \mathcal{C}\}$ in $\mathcal{S}_\mathcal{P}^*$, i.e. $\sigma_{P_\mathcal{C}} \in \mathcal{M}_{P_\mathcal{C}}$ for all $\mathcal{C} \in \mathcal{P}$.

$$\mathcal{R}(s_N) - \mathcal{R}^* = \int_{x \in \mathcal{X}} \left\{ L_{P_x}(s_N(x)) - L_{P_x}^* \right\} \mu(dx) = \sum_{\mathcal{C} \in \mathcal{P}} \int_{x \in \mathcal{C}} \left\{ L_{P_x}(\sigma_{\mathbb{P}_\mathcal{C}}) - L_{P_x}^* \right\} \mu(dx).$$

Now, by virtue of assertion $(i)$ of Lemma 4, we have

$$\mathcal{R}(s_\mathcal{P}) - \mathcal{R}^* \leq 2 \sum_{i<j} \sum_{\mathcal{C} \in \mathcal{P}} \int_{x \in \mathcal{C}} |p_{i,j}(C) - p_{i,j}(x)| \, \mu(dx).$$

Now, observe that, for any $\mathcal{C} \in \mathcal{P}$, all $x \in \mathcal{C}$ and $i < j$, it results from Jensen's inequality and Assumption 3 that

$$|p_{i,j}(\mathcal{C}) - p_{i,j}(x)| \leq \int_{x \in \mathcal{C}} |p_{i,j}(x') - p_{i,j}(x)| \, \mu(dx')/\mu(\mathcal{C}) \leq M \delta_\mathcal{P},$$

which establishes (22).

We now prove the second assertion. For any measurable set $\mathcal{C} \subset \mathcal{X}$ such that $\mu(\mathcal{C}) > 0$, we set $p_{i,j}(\mathcal{C}) = \mathbb{P}\{\Sigma(i) < \Sigma(j) \mid X \in \mathcal{C}\}$ for $i < j$. Suppose that $x \in \mathcal{C}_k$, $k \in \{1, \ldots, K\}$. It follows from assertion (ii) in Lemma 4 combined with Jensen's inequality and Assumption

3 that:

$$d_\tau\left(\sigma^*_{P_x}, s^*_{\mathcal{P}}(x)\right) = d_\tau\left(\sigma^*_{P_x}, \sigma^*_{P_{\mathcal{C}_k}}\right) \leq (1/H)\sum_{i<j}|p_{i,j}(x) - p_{i,j}(\mathcal{C}_k)|$$

$$\leq (1/H)\sum_{i<j}\mathbb{E}\left[|p_{i,j}(x) - p_{i,j}(X)| \mid X \in \mathcal{C}_k\right] \leq (M/H)\sup_{x' \in \mathcal{C}_k}||x - x'|| \leq (M/H)\cdot\delta_{\mathcal{P}}.$$

**Proof of Theorem 12**

We start with proving the first assertion and consider a RMR rule $s_N$ of the form (24). With the notations of Theorem 9, we have the following decomposition:

$$\mathcal{R}(s_N) - \mathcal{R}^* = (\mathcal{R}(s_N) - \mathcal{R}(s_{\mathcal{P}_N})) + (\mathcal{R}(s_{\mathcal{P}_N}) - \mathcal{R}^*). \tag{45}$$

Consider first the second term on the right hand side of the equation above. It results from the argument of Theorem 9's that:

$$\mathcal{R}(s_{\mathcal{P}_N}) - \mathcal{R}^* \leq M\delta_{\mathcal{P}_N} \to 0 \text{ in probability as } N \to \infty. \tag{46}$$

We now turn to the first term. Notice that, by virtue of Lemma 4,

$$\mathcal{R}(s_N) - \mathcal{R}(s_{\mathcal{P}_N}) = \sum_{\mathcal{C} \in \mathcal{P}_N}\left\{L_{P_\mathcal{C}}(\widehat{\sigma}_\mathcal{C}) - L^*_{P_\mathcal{C}}\right\}\mu(\mathcal{C}) \leq 2\sum_{i<j}\sum_{\mathcal{C} \in \mathcal{P}_N}|\widehat{p}_{i,j}(\mathcal{C}) - p_{i,j}(\mathcal{C})|\mu(\mathcal{C}), \tag{47}$$

where, for any $i<j$ and all measurable $\mathcal{C} \subset \mathcal{X}$, we set

$$\widehat{p}_{i,j}(\mathcal{C}) = (1/(N\widehat{\mu}_N(\mathcal{C})))\sum_{k=1}^{N}\mathbb{I}\{X_k \in \mathcal{C}, \Sigma_k(i) < \Sigma_k(j)\}$$

and $\widehat{\mu}_N(\mathcal{C}) = (1/N)\sum_{k=1}^{N}\mathbb{I}\{X_k \in \mathcal{C}\} = N_\mathcal{C}/N$, with the convention that $\widehat{p}_{i,j}(\mathcal{C}) = 0$ when $\widehat{\mu}_N(\mathcal{C}) = 0$. We incidentally point out that the $\widehat{p}_{i,j}(\mathcal{C})$'s are the pairwise probabilities related to the distribution $\widehat{P}_\mathcal{C} = (1/(N\widehat{\mu}_N(\mathcal{C})))\sum_{k:\ X_k\in\mathcal{C}}\delta_{\Sigma_k}$. Observe that for all $i<j$ and $\mathcal{C} \in \mathcal{P}_N$, we have:

$$\mu(\mathcal{C})\left(\widehat{p}_{i,j}(\mathcal{C}) - p_{i,j}(\mathcal{C})\right) = \left\{\frac{1}{N}\sum_{k=1}^{N}\mathbb{I}\{X_k \in \mathcal{C}, \Sigma_k(i) < \Sigma_k(j)\} - \mathbb{E}\left[\mathbb{I}\{X \in \mathcal{C}, \Sigma(i) < \Sigma(j)\}\right]\right\}$$

$$+ \left\{\left(\frac{\mu(\mathcal{C})}{-\widehat{\mu}_N(\mathcal{C}) + \mu(\mathcal{C})} - 1\right)^{-1} \times \frac{1}{N}\sum_{k=1}^{N}\mathbb{I}\{X_k \in \mathcal{C}, \Sigma_k(i) < \Sigma_k(j)\}\right\}.$$

Combining this equality with the previous bound yields

$$\mathcal{R}(s_N) - \mathcal{R}(s_{\mathcal{P}_N}) \leq 2\sum_{i<j}\{A_N(i,j) + B_N/\kappa_N\}, \tag{48}$$

where we set

$$A_N(i,j) = \sup_{\mathcal{P}\in\mathcal{F}_N}\sum_{\mathcal{C}\in\mathcal{P}}\left|\frac{1}{N}\sum_{k=1}^{N}\mathbb{I}\{X_k \in \mathcal{C}, \Sigma_k(i) < \Sigma_k(j)\} - \mathbb{E}\left[\mathbb{I}\{X \in \mathcal{C}, \Sigma(i) < \Sigma(j)\}\right]\right|,$$

$$B_N = \sup_{\mathcal{P}\in\mathcal{F}_N}\sum_{C\in\mathcal{P}}|\widehat{\mu}_N(C) - \mu(C)|$$

The following result is a straightforward application of the VC inequality for data-dependent partitions stated in Theorem 21.1 of Devroye et al. (1996).

**Lemma 17** *Under the hypotheses of Theorem 12, the following bounds hold true:* $\forall \epsilon > 0$, $\forall N \geq 1$,

$$\begin{aligned}
\mathbb{P}\{A_N(i,j) > \epsilon\} &\leq 8\log(\Delta_N(\mathcal{F}_N))e^{-N\epsilon^2/512} + e^{-N\epsilon^2/2}, \\
\mathbb{P}\{B_N > \epsilon\} &\leq 8\log(\Delta_N(\mathcal{F}_N))e^{-N\epsilon^2/512} + e^{-N\epsilon^2/2}.
\end{aligned}$$

The terms $A_N(i,j)$ and $B_N$ are both of order $O_\mathbb{P}(\sqrt{\log(\Delta_N(\mathcal{F}_N))/N})$, as shown by the lemma above. Hence, using Eq. (48) and the assumption that $\kappa_N \to 0$ in probability as $N \to \infty$, so that $1/\kappa_N = o_\mathbb{P}(\sqrt{N/\log \Delta_N(\mathcal{F}_N)})$, we obtain that $\mathcal{R}(s_N) - \mathcal{R}(s_{\mathcal{P}_N}) \to 0$ in probability as $N \to \infty$, which concludes the proof of the first assertion of the theorem.

We now consider the RMR rule (25). Observe that

$$\begin{aligned}
\mathcal{R}(\widetilde{s}_N) - \mathcal{R}(s_{\mathcal{P}_N}) &= \sum_{\mathcal{C} \in \mathcal{P}_N} \left\{ L_{P_\mathcal{C}}(\widetilde{\sigma}^*_{\widehat{P}_\mathcal{C}}) - L^*_{P_\mathcal{C}} \right\} \mu(\mathcal{C}) \\
&= \sum_{\mathcal{C} \in \mathcal{P}_N} \mathbb{I}\{\widehat{P}_\mathcal{C} \in \mathcal{T}\} \left\{ L_{P_\mathcal{C}}(\sigma^*_{\widehat{P}_\mathcal{C}}) - L^*_{P_\mathcal{C}} \right\} \mu(\mathcal{C}) + \sum_{\mathcal{C} \in \mathcal{P}_N} \mathbb{I}\{\widehat{P}_\mathcal{C} \notin \mathcal{T}\} \left\{ L_{P_\mathcal{C}}(\widetilde{\sigma}^*_{\widehat{P}_\mathcal{C}}) - L^*_{P_\mathcal{C}} \right\} \mu(\mathcal{C}) \\
&\leq \mathcal{R}(s_N) - \mathcal{R}(s_{\mathcal{P}_N}) + \frac{n(n-1)}{2} \sum_{\mathcal{C} \in \mathcal{P}_N} \mathbb{I}\{\widehat{P}_\mathcal{C} \notin \mathcal{T}\} \mu(\mathcal{C}). \quad (49)
\end{aligned}$$

Recall that it has been proved previously that $\mathcal{R}(s_N) - \mathcal{R}(s_{\mathcal{P}_N}) \to 0$ in probability as $N \to \infty$. Observe in addition that

$$\mathbb{I}\{\widehat{P}_\mathcal{C} \notin \mathcal{T}\} \leq \mathbb{I}\{P_\mathcal{C} \notin \mathcal{T}\} + \mathbb{I}\{\widehat{P}_\mathcal{C} \notin \mathcal{T} \text{ and } P_\mathcal{C} \in \mathcal{T}\}$$

and, under Assumption 2,

$$\begin{aligned}
\{P_\mathcal{C} \notin \mathcal{T}\} &\subset \{\delta_{\mathcal{P}_N} \geq M/H\}, \\
\{\widehat{P}_\mathcal{C} \notin \mathcal{T} \text{ and } P_\mathcal{C} \in \mathcal{T}\} &\subset \cup_{i<j}\{|\widehat{p}_{i,j}(\mathcal{C}) - p_{i,j}(\mathcal{C})| \geq H\},
\end{aligned}$$

so that $\sum_{\mathcal{C} \in \mathcal{P}_N} \mathbb{I}\{\widehat{P}_\mathcal{C} \notin \mathcal{T}\}\mu(\mathcal{C})$ is bounded by

$$\begin{aligned}
\mathbb{I}\{\delta_{\mathcal{P}_N} \geq M/H\} &+ \sum_{i<j} \sum_{\mathcal{C} \in \mathcal{P}_N} \mathbb{I}\{|\widehat{p}_{i,j}(\mathcal{C}) - p_{i,j}(\mathcal{C})| \geq H\}\mu(\mathcal{C}) \\
&\leq \mathbb{I}\{\delta_{\mathcal{P}_N} \geq M/H\} + \sum_{i<j} \sum_{\mathcal{C} \in \mathcal{P}_N} |\widehat{p}_{i,j}(\mathcal{C}) - p_{i,j}(\mathcal{C})|\mu(\mathcal{C})/H \\
&\leq \mathbb{I}\{\delta_{\mathcal{P}_N} \geq M/H\} + \frac{1}{H} \sum_{i<j} \{A_N(i,j) + B_N/\kappa_N\},
\end{aligned}$$

re-using the argument that previously lead to (48). This bound clearly converges to zero in probability, which implies that $\mathcal{R}(\widetilde{s}_N) - \mathcal{R}(s_{\mathcal{P}_N}) \to 0$ in probability when combined with (49) and concludes the proof of the second assertion of the theorem.

---
**Ranking Median Regression: the $k$-NN Algorithm**

**Inputs.** Training dataset $\mathcal{D}_N = \{(X_1, \Sigma_1), \ldots, (X_N, \Sigma_N)\}$. Norm $||.||$ on the input space $\mathcal{X} \subset \mathbb{R}^d$. Number $k \in \{1, \ldots, N\}$ of neighbours. Query point $x \in \mathcal{X}$

1. (SORT.) Sort the training points by increasing order of distance to $x$:
$$\|X_{(1,N)} - x\| \leq \ldots \leq \|X_{(N,N)} - x\|.$$

2. (ESTIMATION/APPROXIMATION.) Compute the marginal empirical distribution based on the $k$-nearest neighbors in the input space:
$$\widehat{P}(x) = \frac{1}{k} \sum_{l=1}^{k} \delta_{\Sigma_{(k,N)}}$$

**Output.** Compute the local consensus in order to get the prediction at $x$:
$$s_{k,N}(x) = \widetilde{\sigma}^*_{\widehat{P}(x)}.$$

---

Figure 2: Pseudo-code for the $k$-NN algorithm.

**The $k$-NN algorithm for Ranking Median Regression**

**Proof of Theorem 13**

Denote by $\widehat{p}_{i,j}(x)$'s the pairwise probabilities related to distribution $\widehat{P}(x)$. It follows from Lemma 4 combined with Jensen's inequality, that

$$\mathbb{E}\left[\mathcal{R}(s_{k,N}) - \mathcal{R}^*\right] = \mathbb{E}\left[\int_{x \in \mathcal{X}} (L_{P_x}(s_{k,N}(x)) - L^*_{P_x})\mu(dx)\right] \leq 2 \sum_{i<j} \int_{x \in \mathcal{X}} \mathbb{E}\left[|p_{i,j}(x) - \widehat{p}_{i,j}(x)|\right]$$
$$\leq 2 \sum_{i<j} \int_{x \in \mathcal{X}} \left(\mathbb{E}\left[(p_{i,j}(x) - \widehat{p}_{i,j}(x))^2\right]\right)^{1/2}$$

Following the argument of Theorem 6.2's proof in Györfi et al. (2006), write:

$$\mathbb{E}\left[(\widehat{p}_{i,j}(x) - p_{i,j}(x))^2\right] = \mathbb{E}\left[(\widehat{p}_{i,j}(x) - \mathbb{E}\left[\widehat{p}_{i,j}(x)|X_1, \ldots, X_N\right])^2\right]$$
$$+ \mathbb{E}\left[(\mathbb{E}\left[\widehat{p}_{i,j}(x)|X_1, \ldots, X_N\right] - p_{i,j}(x))^2\right] = I_1(x) + I_2(x).$$

The first term can be upper bounded as follows:

$$I_1(x) = \mathbb{E}\left[\left(\frac{1}{k}\sum_{l=1}^{k}\left(\mathbb{I}\left\{\Sigma_{(l,N)}(i) < \Sigma_{(l,N)}(j)\right\} - p_{i,j}(X_{(l,N)})\right)\right)^2\right]$$

$$= \mathbb{E}\left[\frac{1}{k^2}\sum_{l=1}^{k}Var(\mathbb{I}\{\Sigma(i) < \Sigma(j)\}\,|\,X = X_{(l,N)})\right] \leq \frac{1}{4k}.$$

For the second term, we use the following result.

**Lemma 18** *(Lemma 6.4, Györfi et al. (2006)) Assume that the r.v. $X$ is bounded. If $d \geq 3$, then:*

$$\mathbb{E}\left[\|X_{(1,N)}(x) - x\|^2\right] \leq \frac{c_1}{N^{2/d}},$$

*where $c_1$ is a constant that depends on $\mu$'s support only.*

Observe first that, following line by line the argument of Theorem 6.2's proof in Györfi et al. (2006) (see p.95 therein), we have:

$$I_2(x) = \mathbb{E}\left[\frac{1}{k}\left(\sum_{l=1}^{k}\left(p_{i,j}(X_{(l,N)}) - p_{i,j}(x)\right)\right)^2\right] \leq \mathbb{E}\left[\left(\frac{1}{k}\sum_{l=1}^{k}M\|X_{(l,N)} - x\|\right)^2\right]$$

$$\leq M^2 \mathbb{E}\left[\|X_{(1,\lfloor N/k\rfloor)}(x) - x\|^2\right].$$

Next, by virtue of Lemma 18, we have:

$$\frac{1}{M^2}\lfloor N/k\rfloor^{2/d}\int_{x\in\mathcal{X}}I_2(x)\mu(dx) \leq c_1.$$

Finally, we have:

$$\mathbb{E}\left[\mathcal{R}(s_{k,N}) - \mathcal{R}^*\right] \leq 2\sum_{i<j}\int_{x\in\mathcal{X}}\sqrt{I_1(x) + I_2(x)}\mu(dx)$$

$$\leq \frac{n(n-1)}{2}\left(\frac{1}{\sqrt{k}} + 2\sqrt{c_1}M\left(\frac{k}{N}\right)^{1/d}\right).$$

We now consider the problem of bounding the expectation of the excess of risk of the RMR rule $\widetilde{s}_{k,N}$. Observing that $s_{k,N}(x) = \widetilde{s}_{k,N}(x)$ when $\widehat{P}(x) \in \mathcal{T}$, we have:

$$\mathbb{E}\left[\mathcal{R}(\widetilde{s}_{k,N}) - \mathcal{R}^*\right] = \mathbb{E}\left[\int_{x\in\mathcal{X}}\mathbb{I}\{\widehat{P}(x) \in \mathcal{T}\}(L_{P_x}(\widetilde{s}_{k,N}(x)) - L^*_{P_x})\mu(dx)\right] +$$

$$\mathbb{E}\left[\int_{x\in\mathcal{X}}\mathbb{I}\{\widehat{P}(x) \notin \mathcal{T}\}(L_{P_x}(\widetilde{s}_{k,N}(x)) - L^*_{P_x})\mu(dx)\right]$$

$$\leq \mathbb{E}\left[\mathcal{R}(s_{k,N}) - \mathcal{R}^*\right] + \frac{n(n-1)}{2}\mathbb{E}\left[\int_{x\in\mathcal{X}}\mathbb{I}\{\widehat{P}(x) \notin \mathcal{T}\}\mu(dx)\right].$$

Notice in addition that, under Assumption 2, we have, for all $x \in \mathcal{X}$,
$$\{\widehat{P}(x) \notin \mathcal{T}\} \subset \cup_{i<j} \{|\widehat{p}_{i,j}(x) - p_{i,j}(x) \geq H|\}, \tag{50}$$
so that
$$\mathbb{I}\{\widehat{P}(x) \notin \mathcal{T}\} \leq \sum_{i<j} \frac{|\widehat{p}_{i,j}(x) - p_{i,j}(x)|}{H}. \tag{51}$$
Hence, the second assertion finally results directly from the bounds established to prove the first one.

Let $S_{x,\epsilon}$ denote the closed ball centered at $x$ of radius *epsilon* $> 0$. For $d \leq 2$, the rates of convergence hold under the following additional conditions on $\mu$ (see Györfi et al. (2006)): there exists $\epsilon_0 > 0$, a non negative $g$ such that for all $x \in \mathbb{R}^d$ and $0 < \epsilon \leq \epsilon_0$, $\mu(S_{x,\epsilon}) > g(x)\epsilon^d$ and $\int 1/g(x)^{2/d} \mu(dx) < \infty$.

**The CRIT algorithm**

---

<div align="center">THE CRIT ALGORITHM</div>

**Inputs.** Training dataset $\mathcal{D}_N = \{(X_1, \Sigma_1), \ldots, (X_N, \Sigma_N)\}$. Depth $J \geq 0$. Class of admissible subsets $\mathcal{G}$.

1. (INITIALIZATION.) Set $\mathcal{C}_{0,0} = \mathcal{X}$.

2. (ITERATIONS.) For $j = 0, \ldots, J - 1$ and $k = 0, \ldots, 2^j - 1$:

   Solve
   $$\min_{\mathcal{C} \in \mathcal{G}, \, \mathcal{C} \subset \mathcal{C}_{j,k}} \Lambda_{j,k}(\mathcal{C}),$$
   yielding the region $\mathcal{C}_{j+1,2k}$. Then, set $\mathcal{C}_{j+1,2k+1} = \mathcal{C}_{j,k} \setminus \mathcal{C}_{j+1,2k}$.

3. (LOCAL CONSENSUS.) After $2^J$ iterations, for each terminal cell $\mathcal{C}_{J,k}$ with $k \in \{0, \ldots, 2^J - 1\}$, compute the Kemeny median estimate $\sigma^*_{J,k} = \widetilde{\sigma}^*_{\widehat{P}_{\mathcal{C}_{J,k}}}$.

**Outputs.** Compute the piecewise constant ranking median regression rule:
$$s^*_{T_{2^J}}(x) = \sum_{l=0}^{2^J - 1} \sigma^*_{J,l} \cdot \mathbb{I}\{x \in \mathcal{C}_{J,l}\}.$$

---

Figure 3: Pseudo-code for the CRIT algorithm.

From the original tree $T_{2^J}$, one recursively merges children of a same parent node until the root $T_1$ is reached in a bottom up fashion. Precisely, the *weakest link pruning* consists here in sequentially merging the children $\mathcal{C}_{j+1,2l}$ and $\mathcal{C}_{j+1,2l+1}$ producing the smallest

dispersion increase:
$$\widehat{\mu}_N(\mathcal{C}_{j,l})\gamma_{\widehat{P}_{\mathcal{C}_{j,l}}} - \Lambda_{j,l}(\mathcal{C}_{j+1,2l}).$$

One thus obtains a sequence of ranking median regression trees $T_{2^J} \supset T_{2^J-1} \supset \cdots \supset T_1$, the subtree $T_m$ corresponding to a partition with $\#T_m = m$ cells. The final subtree $T$ is selected by minimizing the complexity penalized intra-cell dispersion:

$$\widetilde{\gamma}_T = \widehat{\gamma}_T + \lambda \times \#T, \tag{52}$$

where $\lambda \geq 0$ is a parameter that rules the trade-off between the complexity of the ranking median regression tree, as measured by $\#T$, and intra-cell dispersion. In practice, model selection can be performed by means of common resampling techniques.

**Remark 19** (EARLY STOPPING) *One stops the splitting process if no improvement can be achieved by splitting the current node $\mathcal{C}_{j,l}$, i.e. if $\min_{\mathcal{C} \in \mathcal{G}} \Lambda(\mathcal{C}) = \sum_{1 \leq k < l \leq N} \mathbb{I}\{(X_k, X_l) \in \mathcal{C}_{j,l}^2\} \cdot d_\tau(\Sigma_k, \Sigma_l)$ (one then set $\mathcal{C}_{j+1,2l} = \mathcal{C}_{j,l}$ by convention), or if a minimum node size, specified in advance, is attained.*

**Remark 20** (ON CLASS $\mathcal{G}$) *The choice of class $\mathcal{G}$ involves a trade-off between computational cost and flexibility: a rich class (of controlled complexity though) may permit to capture the conditional variability of $\Sigma$ given $X$ appropriately but might significantly increase the cost of solving (37). Typically, as proposed in Breiman et al. (1984), subsets can be built by means of axis parallel splits, leading to partitions whose cells are finite union of hyperrectangles. This corresponds to the case where $\mathcal{G}$ is stable by intersection, i.e. $\forall (\mathcal{C}, \mathcal{C}') \in \mathcal{G}^2$, $\mathcal{C} \cap \mathcal{C}' \in \mathcal{G}$) and admissible subsets of any $\mathcal{C} \in \mathcal{G}$ are of the form $\mathcal{C} \cap \{X^{(m)} \geq s\}$ or $\mathcal{C} \cap \{X^{(m)} \leq s\}$, where $X^{(m)}$ can be any component of $X$ and $s \in \mathbb{R}$ any threshold value. In this case, the minimization problem can be efficiently solved by means of a double loop (over the $d$ coordinates of the input vector $X$ and over the data lying in the current parent node), see e.g. Breiman et al. (1984).*

**Interpretability and computational feasability.** The fact that the computation of (local) Kemeny medians takes place at the level of terminal nodes of the ranking median regression tree $T$ only makes the CRIT algorithm very attractive from a practical perspective. In addition, it produces predictive rules that can be easily interpreted by means of a binary tree graphic representation and, when implemented with axis parallel splits, provides, as a by-product, indicators quantifying the impact of each input variable. The relative importance of the variable $X^{(m)}$ can be measured by summing the decreases of empirical $\gamma$-dispersion induced by all splits involving it as splitting variable. More generally, the CRIT algorithm inherits the appealing properties of tree induction methods: it easily adapts to categorical predictor variables, training and prediction are fast and it is not affected by monotone transformations of the predictor variables $X^{(m)}$.

**Aggregation.** Just like other tree-based methods, the CRIT algorithm may suffer from instability, meaning that, due to its hierarchical structure, the rules it produces can be much affected by a small change in the training dataset. As proposed in Breiman (1996), **b**oostrap **agg**regat**ing** techniques may remedy to instability of ranking median regression trees. Applied to the CRIT method, bagging consists in generating $B \geq 1$ bootstrap samples by drawing with replacement in the original data sample and running next the learning algorithm from each of these training datasets, yielding $B$ predictive rules $s_1, \ldots, s_B$. For any prediction point $x$, the ensemble of predictions $s_1(x), \ldots, s_B(x)$ are combined in the sense of Kemeny ranking aggregation, so as to produce a consensus $\bar{s}_B(x)$ in $\mathfrak{S}_n$. Observe that a crucial advantage of dealing with piecewise constant ranking rules is that computing a Kemeny median for each new prediction point can be avoided: one may aggregate the ranking rules rather than the rankings in this case, refer to the analysis below for further details. We finally point out that a certain amount of randomization can be incorporated in each bootstrap tree growing procedure, following in the footsteps of the random forest procedure proposed in Breiman (2001), so as to increase flexibility and hopefully improve accuracy.

**Aggregation of piecewise constant predictive ranking rules**

Let $\mathcal{P}$ be a partition of the input space $\mathcal{X}$. By definition, a subpartition of $\mathcal{P}$ is any partition $\mathcal{P}'$ of $\mathcal{X}$ with the property that, for any $\mathcal{C}' \in \mathcal{P}'$, there exists $\mathcal{C} \in \mathcal{P}$ such that $\mathcal{C}' \subset \mathcal{C}$. Given a collection $\mathcal{P}_1, \ldots, \mathcal{P}_B$ of $B \geq 1$ partitions of $\mathcal{X}$, we call the 'largest' subpartition $\bar{\mathcal{P}}_B$ of the $\mathcal{P}_b$'s the partition of $\mathcal{X}$ that is a subpartition of each $\mathcal{P}_b$ and is such that, any subpartition of all the $\mathcal{P}_b$'s is also a subpartition of $\bar{\mathcal{P}}_B$ (notice incidentally that the cells of $\bar{\mathcal{P}}_B$ are of the form $\mathcal{C}_1 \cap \cdots \cap \mathcal{C}_B$, where $(\mathcal{C}_1, \ldots, \mathcal{C}_B) \in \mathcal{P}_1 \times \cdots \times \mathcal{P}_B$). Considering now $B \geq 1$ piecewise constant ranking rules $s_1, \ldots, s_B$ associated with partitions $\mathcal{P}_1, \ldots, \mathcal{P}_B$ respectively, we observe that the $s_b$'s are constant on each cell of $\bar{\mathcal{P}}_B$. One may thus write: $\forall b \in \{1, \ldots, B\}$,

$$s_b(x) = \sum_{\mathcal{C} \in \bar{\mathcal{P}}_B} \sigma_{\mathcal{C},b} \cdot \mathbb{I}\{x \in \mathcal{C}, b\}.$$

For any arbitrary $s(x) = \sum_{\mathcal{C} \in \bar{\mathcal{P}}_B} \sigma_{\mathcal{C}} \mathbb{I}\{x \in \mathcal{C}\}$ in $\mathcal{S}_{\bar{\mathcal{P}}_B}$, we have:

$$\mathbb{E}_{X \sim \mu}\left[\sum_{b=1}^{B} d_\tau\left(s(X), s_b(X)\right)\right] = \sum_{b=1}^{B} \sum_{\mathcal{C} \in \bar{\mathcal{P}}_B} \mu(\mathcal{C}) d_\tau\left(\sigma_{\mathcal{C}}, \sigma_{\mathcal{C},b}\right) = \sum_{\mathcal{C} \in \bar{\mathcal{P}}_B} \sum_{b=1}^{B} \mu(\mathcal{C}) d_\tau\left(\sigma_{\mathcal{C}}, \sigma_{\mathcal{C},b}\right).$$

Hence, the quantity above is minimum when, for each $\mathcal{C}$, the permutation $\sigma_{\mathcal{C}}$ is a Kemeny median of the probability distribution $(1/B) \sum_{b=1}^{B} \delta_{\sigma_{\mathcal{C},b}}$.

**Consistency preservation and aggregation**

We now state and prove a result showing that the (possibly randomized) aggregation procedure previously proposed is theoretically founded. Precisely, mimicking the argument in Biau et al. (2008) for standard regression/classification and borrowing some of their notations, consistency of ranking rules that are obtained through empirical Kemeny aggregation over a profile of consistent *randomized ranking median regression rules* is investigated. Here, a randomized scoring function is of the form $\mathbf{S}_{\mathcal{D}_n}(., Z)$, where $\mathcal{D}_N = \{(X_1, \Sigma_1), \ldots, (X_N, \Sigma_N)\}$ is the training dataset and $Z$ is a r.v. taking its values in some measurable space $\mathcal{Z}$ that describes the randomization mechanism. A randomized ranking rule $\mathbf{S}_{\mathcal{D}_n}(., Z) : \mathcal{X} \to \mathfrak{S}_N$ is given and consider its ranking median regression risk,

$\mathcal{R}(\mathbf{S}_{\mathcal{D}_N}(., Z))$ namely, which is given by:

$$\mathcal{R}(\mathbf{S}_{\mathcal{D}_N}(., Z)) = \sum_{i<j} \mathbb{P}\left\{ (\mathbf{S}_{\mathcal{D}_N}(X, Z)(j) - \mathbf{S}_{\mathcal{D}_N}(X, Z)(i))(\Sigma(j) - \Sigma(i)) < 0 \right\},$$

where the conditional probabilities above are taken over a random pair $(X, \Sigma)$, independent from $\mathcal{D}_N$. It is said to be consistent iff, as $N \to \infty$,

$$\mathcal{R}(\mathbf{S}_{\mathcal{D}_N}(., Z)) \to \mathcal{R}^*,$$

in probability. When the convergence holds with probability one, one says that the randomized ranking rule is strongly consistent. Fix $B \geq 1$. Given $\mathcal{D}_N$, one can draw $B$ independent copies $Z_1, \ldots, Z_B$ of $Z$, yielding the ranking rules $\mathbf{S}_{\mathcal{D}_N}(., Z_b)$, $1 \leq b \leq B$. Suppose that the ranking rule $\bar{\mathbf{S}}_B(.)$ minimizes

$$\sum_{b=1}^{B} \mathbb{E}\left[ d_\tau\left(s(X), \mathbf{S}_{\mathcal{D}_N}(X, Z_b)\right) \mid \mathcal{D}_N, Z_1, \ldots, Z_B \right] \tag{53}$$

over $s \in \mathcal{S}$. The next result shows that, under Assumption 2 (strong) consistency is preserved for the ranking rule $\bar{\mathbf{S}}_B(X)$ (and the convergence rate as well, by examining its proof).

**Theorem 21** (CONSISTENCY AND AGGREGATION.) *Assume that the randomized ranking rule $\mathbf{S}_{\mathcal{D}_N}(., Z)$ is consistent (respectively, strongly consistent) and suppose that Assumption 2 is satisfied. Let $B \geq 1$ and, for all $N \geq 1$, let $\bar{\mathbf{S}}_B(x)$ be a Kemeny median of $B$ independent replications of $\mathbf{S}_{\mathcal{D}_N}(x, Z)$ given $\mathcal{D}_N$. Then, the aggregated ranking rule $\bar{\mathbf{S}}_B(X)$ is consistent (respectively, strongly consistent).*

**Proof** Let $s^*(x)$ be an optimal ranking rule, i.e $\mathcal{R}(s^*) = \mathcal{R}^*$. Observe that, with probability one, we have: $\forall b \in \{1, \ldots, B\}$,

$$d_\tau\left(s^*(X), \bar{\mathbf{S}}_B(X)\right) \leq d_\tau\left(s^*(X), \mathbf{S}_{\mathcal{D}_N}(X, Z_b)\right) + d_\tau\left(\mathbf{S}_{\mathcal{D}_N}(X, Z_b), \bar{\mathbf{S}}_B(X)\right),$$

and thus, by averaging,

$$d_\tau\left(s^*(X), \bar{\mathbf{S}}_B(X)\right) \leq \frac{1}{B}\sum_{b=1}^{B} d_\tau\left(s^*(X), \mathbf{S}_{\mathcal{D}_N}(X, Z_b)\right) + \frac{1}{B}\sum_{b=1}^{B} d_\tau\left(\mathbf{S}_{\mathcal{D}_N}(X, Z_b), \bar{\mathbf{S}}_B(X)\right)$$

$$\leq \frac{2}{B}\sum_{b=1}^{B} d_\tau\left(s^*(X), \mathbf{S}_{\mathcal{D}_N}(X, Z_b)\right).$$

Taking the expectation w.r.t. to $X$, one gets:

$$\mathcal{R}(\bar{\mathbf{S}}_B) - \mathcal{R}^* \leq \mathbb{E}\left[ d_\tau\left(s^*(X), \bar{\mathbf{S}}_B(X)\right) \mid \mathcal{D}_N, Z_1, \ldots, Z_B \right]$$
$$\leq \frac{1}{B}\sum_{b=1}^{B} \mathbb{E}\left[ d_\tau\left(s^*(X), \mathbf{S}_{\mathcal{D}_N}(X, Z_b)\right) \mid \mathcal{D}_N, Z_1, \ldots, Z_B \right] \leq \frac{1}{BH}\sum_{b=1}^{B} \left\{ \mathcal{R}(\mathbf{S}_{\mathcal{D}_N}(., Z_b)) - \mathcal{R}^* \right\}$$

using (13) and Assumption 2. This bound combined with the (strong) consistency assumption yields the desired result. ∎

Of course, the quantity (53) is unknown in practice, just like distribution $\mu$, and should be replaced by a statistical version based on an extra sample composed of independent copies of $X$. Statistical guarantees for the minimizer of the empirical variant over a subclass $\mathcal{S}_0 \subset \mathcal{S}$ of controlled complexity (*i.e.* under Assumption 1) could be established by adapting the argument of Theorem 21's proof. However, from a practical perspective, aggregating the predictions (at a given point $x \in \mathcal{X}$) rather than the predictors would be much more tractable from a computational perspective (since it would permit to exploit the machinery for approximate Kemeny aggregation of subsection 2.2). This will be investigated in a forthcoming paper.

| Dataset distribution | Setting 1 | | | Setting 2 | | | Setting 3 | | |
|---|---|---|---|---|---|---|---|---|---|
| | n=3 | n=5 | n=8 | n=3 | n=5 | n=8 | n=3 | n=5 | n=8 |
| Piecewise constant | 0.0698* | 0.1290* | 0.2670* | 0.0173* | 0.0405* | 0.110* | 0.0112* | 0.0372* | 0.0862* |
| | 0.0473** | 0.136** | 0.324** | 0.0568** | 0.145** | 0.2695** | 0.099** | 0.1331** | 0.2188** |
| | (0.578) | (1.147) | (2.347) | (0.596) | (1.475) | (3.223) | (0.5012) | (1.104) | (2.332) |
| Mallows with $\phi=2$ | 0.3475 * | 0.569* | 0.9405 * | 0.306* | 0.494* | 0.784* | 0.289* | 0.457* | 0.668* |
| | 0.307** | 0.529** | 0.921** | 0.308** | 0.536** | 0.862** | 0.3374** | 0.5714** | 0.8544** |
| | (0.719) | (1.349) | (2.606) | (0.727) | (1.634) | (3.424) | (0.5254) | (1.138) | (2.287) |
| Mallows with $\phi=1$ | 0.8656* | 1.522* | 2.503* | 0.8305 * | 1.447 * | 2.359* | 0.8105* | 1.437* | 2.189* |
| | 0.7228** | 1.322** | 2.226** | 0.723** | 1.3305** | 2.163** | 0.7312** | 1.3237** | 2.252** |
| | (0.981) | (1.865) | (3.443) | (1.014) | (2.0945) | (4.086) | (0.8504) | (1.709) | (3.005) |

Table 1: Empirical risk averaged on 50 trials on simulated data.